\documentclass{article}
\usepackage{amsfonts}
\usepackage{amsmath}
\usepackage{hyperref}
\usepackage{curves}
\usepackage{enumerate}
\usepackage[utf8]{inputenc}
\usepackage{latexsym}
\usepackage{braket}
\usepackage{mathrsfs}
\usepackage{mathtools}
\usepackage{color}
\usepackage{xcolor}
\usepackage{graphicx}
\usepackage{subcaption}
\usepackage{cleveref}
\usepackage{enumitem}
\usepackage[toc,page]{appendix}
\usepackage[a4paper, total={6in, 8in}]{geometry}

\usepackage{cite}
\usepackage[english]{babel}
\usepackage{babelbib}
\usepackage{multirow}
\usepackage{tabularx}
\usepackage{makecell}
\setcellgapes{4pt}
\usepackage[thinlines]{easytable}

\usepackage{colortbl}
\def\bal#1\eal{\begin{align}#1\end{align}}
\usepackage{color,soul}
\newcommand{\bsub}{\begin{subequations}}
\newcommand{\esub}{\end{subequations}}
\def\bal#1\eal{\begin{align}#1\end{align}}

\newcommand{\gra}{{\alpha}} \newcommand{\grb}{{\beta}}  
  \newcommand{\grh}{{\eta}} 
 \newcommand{\grk}{{\kappa}} \newcommand{\grl}{{\lambda}} \newcommand{\grm}{{\mu}}
\newcommand{\grn}{{\nu}} \newcommand{\grj}{{\xi}}  
 \newcommand{\grs}{{\sigma}} \newcommand{\grt}{{\tau}} \newcommand{\gry}{{\upsilon}}
   \newcommand{\grv}{{\omega}}

 \newcommand{\grK}{{\rm K}} \newcommand{\grL}{{\Lambda}} \newcommand{\grM}{{\rm M}}
\newcommand{\grN}{{\rm N}}   
\newcommand{\grR}{{\rm R}} \newcommand{\grS}{{\Sigma}}

\begin{document}

\title{\textbf{On solvable Lie algebras and integration method of ordinary differential equations}}
\author{ \textbf{T. Pailas}$^a$\thanks{teopailas879@hotmail.com} ,\,\,\textbf{P. A. Terzis}$^a$\thanks{pterzis@phys.uoa.gr} ,\,\,\textbf{T. Christodoulakis}$^a$\thanks{tchris@phys.uoa.gr}\\
\normalsize
{$^a$\it Department of Nuclear and Particle Physics, Faculty of Physics,}\\
\normalsize{\it National and Kapodistrian University of Athens, Athens 15784, Greece}}

\maketitle

\abstract{Some more general ``inheritance conditions" have been found for a given set of symmetry generators $\{\mathbf{Z}_{\bar{l}}\}$ acting on some set of coupled ordinary differential equations, once the ``first integration method'' has been applied upon some Abelian sub-algebra of generators $\{\mathbf{Z}_{i}\}$ for the system of equations to be reduced. We have proven the following theorem: all the generators of some solvable sub-algebra can be used via the ``first integration method'' in order to reduce (the order) and/or (possibly) integrate the system of equations. The specific order in which the step by step reduction process has to be performed, and the condition to obtain the analytic solution solely in terms of quadratures, are also provided. We define the notion ``optimum'' for some given maximal, solvable sub-algebra and prove a theorem stating that: this algebra leads to the most profitable way of using the symmetry generators for the integration procedure at hand.}

\newpage

\section{Introduction}

The concept of symmetry possesses a fundamental role in the field of mathematical physics. An extensive analysis concerning with the symmetries of fundamental physics can be found in \cite{Sundermeyer:2014kha}. There are many branches corresponding to the different kinds of symmetries. One of these is mostly referred to as, continuous symmetries of differential equations (both ordinary and partial), and is concerned with continuous (and invertible) transformations who are mapping solutions into solutions of the equations at hand. The notion of continuous groups of symmetry transformations (now called Lie groups) was first introduced in a series of seminal works \cite{Lie1,Lie2,Lie1888,Lie1891,Lie5,Lie6}, by Sophus Lie; the purpose was to unify and extend various methods related to topics such as, homogeneous and separable equations, Laplace transformation, reduction of order and more. These symmetry transformations are also called Lie point symmetries due to the mix of independent and dependent variables. 

We will refer only to a fraction of the many applications that can be found in the literature. Construction of special solutions of partial differential equations based on the symmetry transformations method was first noticed by S. Lie himself, but was not studied in detail until the work of Ovsiannikov and his group \cite{Ovs58,Ovs59,Ovs60,Ovs61}. A special kind of solutions referred to as similarity solutions, were studied for some specific partial differential equations, as well as the conditions under which a non-linear differential equation is equivalent to a linear one, in
\cite{B0,B1,B2,B3,B4,B5,Doyle}. E. Noether studied in \cite{Noether1918} the symmetries of action integrals and their relation to conservation laws for Euler-Lagrange equations. These symmetries can be determined through Lie's method since they are also symmetries of the respective equations. 

There are also many applications of Lie's theory in the broad area of general relativity. Lie point symmetries of geodesic equations and collineations have been studied in \cite{Tsamparlis:2011wf} and applied to the following cases: Einstein spaces, Schwarzschild and FLRW spacetimes. A dust fluid arising from solutions of the Wheeler-DeWitt equation, generated by Lie point symmetries in scalar field cosmology, have been found in \cite{Paliathanasis:2017ocj}. Constraints and analytical solutions of $f(R)$ theories of gravity based on Lie point symmetries have been obtained in \cite{Paliathanasis:2011jq}. The solution space of the four dimensional Bianchi Type $I-VII$ spacetimes in vacuum, as well as for the case of five dimensional Bianchi Type I, has been found with the aid of the Automorphisms group. This group constitutes a subset of Lie point symmetries of the Einstein's equations for the Bianchi Type geometry \cite{Christodoulakis:2006vi,Terzis:2008ev,Terzis:2010dk,Pailas:2018tzy}. The Lie point symmetries and the variational symmetries for some cases of minisuperspace Einstein's gravity were studied in \cite{Christodoulakis:2013xha}. 

A generalization of Lie point symmetry transformations has been achieved by A. V. Backlund in his seminal papers \cite{Bac1,Bac2,Bac3,Bac4}. These are transformations mixing the independent, dependent and the derivatives of the dependent variables. A particular case are contact transformations (in which the new independent variables are linear in the derivatives of the old dependent variables), \cite{LE} and some applications of them can be found in \cite{Dimakis:2015rba}. Robert Geroch in \cite{G1,G2} provided a method for generating new solutions of Einstein's field equations from already known ones. The transformation that relates the old and new metrics is a Lie-Backlund transformation, even though he did not follow this route to provide his results. Lie-Backlund transformations are discussed in detail by Anderson and Ibragimov \cite{Ibra1}, Ibragimov \cite{Ibra2} and Olver \cite{olver2000applications}. As was shown by Olver \cite{Olv1}, the invariance of partial differential equations under a Lie-Backlund symmetry, usually leads to the existence of infinite number of such symmetries connected via recursion operators.

There exist also generalizations of Lie's methods for obtaining solutions of the differential equations, in the cases where the equations at hand are difference equations \cite{Mae}. Furthermore, the case of differential-difference equations and the application to the Toda lattice was studied by D. Levi and P. Winterwitz \cite{LEVI1991335}. An algorithm for obtaining the Lie point symmetries of differential equations on fixed non-transforming lattices is presented in \cite{Levi_2010}. Finally, a method is given to derive the point symmetries of partial differential-difference equations in \cite{QUISPEL1992379}. It has been applied to the Kacvan-Moerbeke equation and has been found that the symmetries form a Kac-Moody-Virasoro algebra.

One of the primary and most important applications of Lie point symmetries, is their use for the reduction of order and hence the (partial or total) integration of the equations at hand. There are mainly three different methods to do such a reduction. Some of them are listed and thoroughly examined in \cite{olver2000applications}, others in \cite{stephani_1990} and \cite{bluman2013symmetries}. A prominent problem that all three methods have to deal with, is to provide an answer to the question of how the symmetry generators can be used in the most ``profitable'' way; where ``profitable'' is tantamount to as many generators as possible.

The first method to appear in all three previous books is called ``Normal form of generators in the space of variables''. For simplicity, we will refer to it as the ``first integration method''. The procedure is presented in some detail in section II.

The second method outlined by Stephani \cite{stephani_1990} is called the ``Normal form of generators in the space of first integrals''. As it was there proven, this method can be used when the number of symmetries is higher or equal to the order of the equation. An important, worth mentioning aspect, is that it is necessary to use the generators for the integration procedure in a very specific order. This depends on the properties of the algebra. When there is a solvable sub-algebra of dimension equal to the order of the equation, and furthermore integration is performed in the correct order, the solution to the equation is given solely in terms of quadratures. We note that this is a combination of a specific method of integration, with one of the properties of the algebra at hand, namely solvability. 

The final method is sketched in \cite{stephani_1990} (papes $87-88$), while it is extensively studied in \cite{olver2000applications}. It is referred to as the method of ``Differential invariants''. As it was proven, all the differential invariants of higher order are constructed as derivatives of the first two (zeroth and first order respectively), by considering one as the independent and the other as the dependent variable. This reduction procedure can again be combined with the solvability property of some sub-algebra and a proper order is presented to the purpose of using all the generators of the solvable sub-algebra. Once more, when the dimension of the solvable sub-algebra is equal to the order of the equation at hand, the solution can, in principle, be given in terms of quadratures.  

The drawback of the second method is that it requires the number of symmetries be at least equal to the order of the equation. In many occasions, especially when there is a system of equations instead of only one, this is not the case. The third method, on the other hand, does not presupposes any lower bound for the number of symmetry generators; but it lacks of applicability in the case of system of equations. This is due to the fact that, constructing higher order differential invariants from the first two, is far from trivial and also ambiguous.

The simplest appears to be the ``first integration method''; it can be applied with ease, either when the symmetry generators are less than the order of the equation or when there is some system of equations instead of only one. The results contained in this work are briefly outlined as follows:
\begin{enumerate} 
\item Generalization of the ``inheritance conditions" in the case of existence of more than one generators forming an Abelian sub-algebra. 
\item Combination of the property of solvability and the ``first integration method'', as well as presentation of the proper order in which the generators have to be used. 
\item Finally, providing an answer to the question of whether or not the ``optimum'' maximal solvable sub-algebra is the most profitable way to use the symmetry generators of the entire algebra (based on the ``first integration method'').
\end{enumerate}
To the best of our knowledge, these points have not been previously proven or even discussed in the literature.

The paper is organized as follows: In section II, the ``first integration method'' is presented in detail. Furthermore, the proof for the new ``inheritance conditions" is given along with other side effects. All these results are gathered in Theorem 1 and Corollary 1. The section III is devoted to recalling the definition of solvable algebras, alongside with the ones of derived and ``coset" algebras. The combination of solvability and the ``first integration method'' is presented in section IV, where the proper order of reduction is given in Theorem 5. An example is provided in section V. The proof that the maximal solvable sub-algebra is the most profitable way to use the symmetry generators is dispensed in section VI. Alongside, a definition of the proper basis of the generators is given, in which we acquire the ``optimum'' (maximal dimensional) solvable sub-algebra;. A discussion can be found in section VII and finally an Appendix is included.

\section{``First integration method''}

All the statements below, can be easily generalized for systems of higher order equations. For simplicity, we will consider a system of second order ordinary differential equations written in the form
\begin{align}
\ddot{y}^{\grm}=\grv^{\grm}(x,y^{\grn},\dot{y}^{\grn}),\label{e1}
\end{align}
where $x$ is the common independent variable, $y^{\grn}$ the dependent variables, $\dot{y}^{\grn}=\frac{dy^{\grn}}{dx}\,,\ddot{y}^{\grn}=\frac{d\dot{y}^{\grn}}{dx}$, and the indices $\grm,\grn, \grk$ run from 1 to N.
As can be found in \cite{stephani_1990}, considering first integrals $f^{\grm}(x,y^{\grn},\dot{y}^{\grn})$ of \eqref{e1} and taking
\begin{align}
\frac{df^{\grm}(x,y^{\grn},\dot{y}^{\grn})}{dx}=0,\label{en1}
\end{align}  
one can extract the operator $\mathbf{A}$
\begin{align}
\mathbf{A}=\partial_{x}+\dot{y}^{\grk}\partial_{y^{\grk}}+\grv^{\grk}(x,y^{\grn},\dot{y}^{\grn})\partial_{\dot{y}^{\grk}},\label{e2}
\end{align} 
so as to represent the system of equations as
\begin{align}
\mathbf{A}f^{\grm}(x,y^{\grn},\dot{y}^{\grn})=0.\label{e3565}
\end{align}
At this point we assume the existence of a number $r$ of Lie point symmetries represented by the following linearly independent, infinitesimal generators
\begin{align}
\mathbf{Z}_{\grM}=\grj_{\grM}(x,y^{\grn})\partial_{x}+\grh^{\grm}_{\grM}(x,y^{\grn})\partial_{y^{\grm}},\label{e3}
\end{align}
for some arbitrary but given set of functions $\grj_{\grM}(x,y^{\grn}),\,\grh^{\grm}_{\grM}(x,y^{\grn})$, where the capital $M,N,K$ acquire the values $M,N,K=1,...r$. Given the generators, their first prolongation(we only need the first prolongation for our purposes) can be calculated via the following formula
\begin{align}
\mathbf{Z}^{(1)}_{\grM}=\grj_{\grM}(x,y^{\grn})\partial_{x}+\grh^{\grm}_{\grM}(x,y^{\grn})\partial_{y^{\grm}}+\grh^{(1)\grm}_{\grM}(x,y^{\grn},\dot{y}^{\grn})\partial_{\dot{y}^{\grm}},\label{e4}
\end{align}
where 
\begin{align}
\grh^{(1)\grm}_{\grM}(x,y^{\grn},\dot{y}^{\grn})=\frac{d \grh^{\grm}_{\grM}(x,y^{\grn})}{dx}-\dot{y}^{\grm}\frac{d\grj_{M}(x,y^{\grn})}{dx}\label{e5}.
\end{align}
The symmetry generators form a Lie algebra 
\begin{align}
\left[\mathbf{Z}_{\grM},\mathbf{Z}_{\grN}\right]=C^{\grL}_{\grM\grN}\mathbf{Z}_{\grL}\label{e6},
\end{align} 
with $C^{\grL}_{\grM\grN}$ the structure constants of the algebra which satisfy the following relations
\begin{enumerate}
\item $C^{\grL}_{\grM\grN}=-C^{\grL}_{\grN\grM}$,
\item $C^{\grK}_{\grM\grN}C^{\grR}_{\grS\grK}+C^{\grK}_{\grN\grS}C^{\grR}_{\grM\grK}+C^{\grK}_{\grS\grM}C^{\grR}_{\grN\grK}=0$,
\item $C^{\grK}_{\grN\grS}C^{\grM}_{\grM\grK}=0,$
\end{enumerate}
and $\left[\cdot,\cdot\right]$ the commutator. The structure constants are the same for the prolongation of the generators as well.
 
The second formulation of the symmetry condition for a system of second order ordinary differential equations reads, according to \cite{stephani_1990}
\begin{align}
\left[\mathbf{Z}^{(1)}_{\grM},\mathbf{A}\right]=-\left(\mathbf{A}(\grj_{\grM})\right)\mathbf{A}\equiv\grl_{\grM}\mathbf{A}.\label{e7}
\end{align}
Based on the ``first integration method'' (i.e. transform the generators into their normal form in the space of variables) \cite{stephani_1990}, to reduce the order for one of the differential equations \eqref{e1}, we have to transform one of the symmetry generators, let us say $\mathbf{Z}_{1}$, into it's normal form. That is equivalent to search for independent and dependent coordinates $t(x,y^{\grn})$, $s^{\grm}(x,y^{\grn})$ respectively, such that in these coordinates
\begin{align}
\mathbf{Z}_{1}=\partial_{s^{1}}.\label{e8}
\end{align} 
When this is achieved, the system of the transformed equations does not depend on $s^{1}$ and thus the order is reduced via the restriction to the hypersurfaces $s^{1}$=constant, $w^{1}=\dot{s}^{1}$, $\dot{w}^{1}=\ddot{s}^{1}$. In general, we are not sure which of the rest of the symmetries will be inherited in the reduced equations. As it was proven in \cite{stephani_1990}, the generators $\mathbf{Z}_{\grM}$ of those symmetries must satisfy the relation
\begin{align}
\left[\mathbf{Z}_{1},\mathbf{Z}_{\grM}\right]=C^{1}_{1\grM}\mathbf{Z}_{1},\label{e9}
\end{align}
where $M=\grs,...,r$ and $\grs\geq2$. Thus, it is to our best interest to choose that $\mathbf{Z}_{1}$ which satisfies the relation \eqref{e9} with the maximal number of $\mathbf{Z}_{\grM}$ so as to inherit the maximal symmetry group.

In this work we first generalize the condition \eqref{e9} in the case where a number of generators form an Abelian sub-algebra, where we can transform all it's generators into normal form at once.\\ 
The following theorem holds:

\textbf{\underline{Theorem 1:}}
\textbf{\textit{Given a system of ordinary differential equations represented in terms of a partial differential operator $\mathbf{A}$ of first order, we assume the existence of r number of linearly independent Lie point symmetries with generators $\mathbf{Z}_{M}$, $M=1,...,r$. Consider also A to be linearly independent of those generators. Furthermore, assume there exists an $m\leq r$ dimensional Abelian sub-algebra with generators $\mathbf{Z}_{i}$, which acts transitive in a subspace of dimension m. When these generators are transformed into their normal form and thus the system is reduced, the generators $\mathbf{Z}_{\bar{l}}$, to be inherited are those satisfying the following conditions:}}
\begin{enumerate}
\item 
$\left[\mathbf{Z}_{i},\mathbf{Z}_{\bar{l}}\right]=C_{i\bar{l}}^{l}\mathbf{Z}_{l}$.\\
\textbf{\textit{ Furthermore, the reduced objects}} $\hat{\mathbf{A}},\,\mathbf{Y}_{\bar{l}}$ \textbf{\textit{in the hypersurfaces $s^{j}$=constant, $u^{j}=\dot{s}^{j}$, will satisfy the relations below}}
\item $\mathbf{\hat{A}}=\partial_{t}+\dot{s}^{\grb}\partial_{s^{\grb}}+\hat{\grv}^{i}(t,s^{\gra},u^{j},\dot{s}^{\gra})\partial_{u^{i}}+\hat{\grv}^{\grb}(t,s^{\gra},u^{j},\dot{s}^{\gra})\partial_{\dot{s}^{\grb}}$, 
\item $\mathbf{Y}_{\bar{l}}=\grj_{\bar{l}}(t,s^{\gra})\partial_{t}+\grh_{\bar{l}}^{\grb}(t,s^{\gra})\partial_{s^{\grb}}+\grh^{(1)i}_{\bar{l}}(t,s^{\gra},u^{j},\dot{s}^{\gra})\partial_{u^{i}}+\grh^{(1)\grb}_{\bar{l}}(t,s^{\gra},u^{j},\dot{s}^{\gra})\partial_{\dot{s}^{\grb}}$,
\item $[\mathbf{Y}_{\bar{l}},\mathbf{Y}_{\bar{k}}]=C^{\bar{q}}_{\bar{l}\bar{k}}\mathbf{Y}_{\bar{q}}$,
\end{enumerate}

\textbf{\underline{Proof}}

Since the $m$-dimensional sub-algebra, $m\leq r$, is assumed to be Abelian  the generators $\mathbf{Z}_{i}$ satisfy
\begin{align}
\left[\mathbf{Z}_{i},\mathbf{Z}_{j}\right]=0\label{e11}.
\end{align}
Accordingly, the assumption of transitivity in a subspace of dimension m (in Appendix A we explain why we need the transitivity), implies that there exists independent and dependent coordinates $t(x,y^{\grn})$, $s^{\grm}(x,y^{\grn})$ respectively, such that 
\begin{align}
\mathbf{Z}_{i}=\partial_{s^{i}},\label{e12}
\end{align}
while the equations \eqref{e1} are transformed into the set
\begin{align}
&\ddot{s}^{i}=\tilde{\grv}^{i}(t,s^{j},s^{\gra},\dot{s}^{j},\dot{s}^{\gra}),\label{e13}\\
&\ddot{s}^{\grb}=\tilde{\grv}^{\grb}(t,s^{j},s^{\gra},\dot{s}^{j},\dot{s}^{\gra}),\label{e14}
\end{align} 
where $i,j=1,...,m$ and $\gra,\grb=m+1,...,N$, where a reshuffling of the indices may be employed. The partial differential operator of first order corresponding to the equations \eqref{e13},\eqref{e14} as well as the first prolongation of the symmetry generators read
\begin{align}
&\mathbf{A}=\partial_{t}+\dot{s}^{i}\partial_{s^{i}}+\dot{s}^{\grb}\partial_{s^{\grb}}+\tilde{\grv}^{i}(t,s^{j},s^{\gra},\dot{s}^{j},\dot{s}^{\gra})\partial_{\dot{s}^{i}}+\tilde{\grv}^{\grb}(t,s^{j},s^{\gra},\dot{s}^{j},\dot{s}^{\gra})\partial_{\dot{s}^{\grb}},\label{e15}\\
&\mathbf{Z}^{(1)}_{l}=\partial_{s^{l}},\label{e16}
\end{align}
\begin{align}
\mathbf{Z}^{(1)}_{\bar{l}}=&\grj_{\bar{l}}(t,s^{j},s^{\gra})\partial_{t}+\grh^{i}_{\bar{l}}(t,s^{j},s^{\gra})\partial_{s^{i}}+\grh^{\grb}_{\bar{l}}(t,s^{j},s^{\gra})\partial_{s^{\grb}}+\grh^{(1)i}_{\bar{l}}(t,s^{j}\,s^{\gra},\dot{s}^{j},\dot{s}^{\gra})\partial_{\dot{s}^{i}}\nonumber\\
&+\grh^{(1)\grb}_{\bar{l}}(t,s^{j}\,s^{\gra},\dot{s}^{j},\dot{s}^{\gra})\partial_{\dot{s}^{\grb}},\label{e17}
\end{align}
where $l=1,...,m$ and $\bar{l}=m+1,...,\tilde{m}$. The equations \cref{e13,e14,e15,e16,e17} are just the expressions of the  transformed objects under the coordinate transformation that we performed. The symmetry condition for $\mathbf{Z}^{(1)}_{l}$ implies
\begin{align}
&\left[\mathbf{Z}^{(1)}_{l},\mathbf{A}\right]=-\underbrace{\left(\mathbf{A}(\grj_{l})\right)}_{=0}\mathbf{A}\Rightarrow\nonumber\\
&\mathbf{Z}^{(1)}_{l}\left(\mathbf{A}\right)-\mathbf{A}\left(\mathbf{Z}^{(1)}_{l}\right)=0\Rightarrow\nonumber\\
&\partial_{s^{l}}\left[\tilde{\grv}^{i}(t,s^{j},s^{\gra},\dot{s}^{j},\dot{s}^{\gra})\right]\partial_{\dot{s}^{i}}+\partial_{s^{l}}\left[\tilde{\grv}^{\grb}(t,s^{j},s^{\gra},\dot{s}^{j},\dot{s}^{\gra})\right]\partial_{\dot{s}^{\grb}}=0\Rightarrow\nonumber
\end{align}
\begin{align}
&\tilde{\grv}^{i}(t,s^{j},s^{\gra},\dot{s}^{j},\dot{s}^{\gra})=\hat{\grv}^{i}(t,s^{\gra},\dot{s}^{j},\dot{s}^{\gra}),\label{e18}\\
&\tilde{\grv}^{\grb}(t,s^{j},s^{\gra},\dot{s}^{j},\dot{s}^{\gra})=\hat{\grv}^{\grb}(t,s^{\gra},\dot{s}^{j},\dot{s}^{\gra}),\label{e19}
\end{align}
where the independence of the basis vectors $\partial_{\dot{s}^{i}}$, $\partial_{\dot{s}^{\grb}}$ has been used and some intermediate calculations between the second and third line has been omitted in the interest of brevity. Thus, the equations \eqref{e13},\eqref{e14} are actually
\begin{align}
&\ddot{s}^{i}=\hat{\grv}^{i}(t,s^{\gra},\dot{s}^{j},\dot{s}^{\gra})\label{e20},\\
&\ddot{s}^{\grb}=\hat{\grv}^{\grb}(t,s^{\gra},\dot{s}^{j},\dot{s}^{\gra})\label{e21}.
\end{align}
By the restriction to the hypersurfaces $s^{i}$=constant, $u^{i}=\dot{s}^{i}$, the system of the original second order ordinary differential equations \eqref{e1} has been reduced to a system of first order differential equations coupled to a system of second order ordinary differential equations and a set of quadratures.
\begin{align}
&\dot{u}^{i}=\hat{\grv}^{i}(t,s^{\gra},u^{j},\dot{s}^{\gra}),\label{e22}\\
&\ddot{s}^{\grb}=\hat{\grv}^{\grb}(t,s^{\gra},u^{j},\dot{s}^{\gra}),\label{23}\\
&s^{i}=\int{u^{i}dt}+c^{i},\label{e24}
\end{align}
where $c^{i}$ are some integration constants. In order for the procedure of reduction to proceed, we introduce the coordinates attached to the hypersurfaces $s^{i}=constant$ to be $(t,s^{\grb},u^{i},\dot{s}^{\grb})$.  The partial differential operator and the first prolongation of the rest of the symmetries have to be calculated in the reduced coordinates. On way to do so is the following:
\begin{align}
&\hat{\mathbf{A}}=\mathbf{A}(t)\partial_{t}+\mathbf{A}(u^{i})\partial_{u^{i}}+\mathbf{A}(s^{\grb})\partial_{s^{\grb}}+\mathbf{A}(\dot{u}^{i})\partial_{\dot{u}^{i}}+\mathbf{A}(\dot{s}^{\grb})\partial_{\dot{s}^{\grb}}\label{e25},\\
&\mathbf{Y}_{\bar{l}}=\mathbf{Z}^{(1)}_{\bar{l}}(t)\partial_{t}+\mathbf{Z}^{(1)}_{\bar{l}}(u^{i})\partial_{u^{i}}+\mathbf{Z}^{(1)}_{\bar{l}}(s^{\grb})\partial_{s^{\grb}}+\mathbf{Z}^{(1)}_{\bar{l}}(\dot{u}^{i})\partial_{\dot{u}^{i}}+\mathbf{Z}^{(1)}_{\bar{l}}(\dot{s}^{\grb})\partial_{\dot{s}^{\grb}},\label{e26}
\end{align}
where we have to use the expressions of the coordinates to the hypersurfaces in terms of the original coordinates. After some trivial mathematical calculations we end up with the expressions
\begin{align}
&\hat{\mathbf{A}}=\partial_{t}+\dot{s}^{\grb}\partial_{s^{\grb}}+\hat{\grv}^{i}(t,\dot{s}^{j},s^{\gra},\dot{s}^{\gra})\partial_{\dot{s}^{i}}+\hat{\grv}^{\grb}(t,\dot{s}^{j},s^{\gra},\dot{s}^{\gra})\partial_{\dot{s}^{\grb}},\label{e27}\\
&\mathbf{Y}_{\bar{l}}=\grj_{\bar{l}}(t,s^{j},s^{\gra})\partial_{t}+\grh^{\grb}_{\bar{l}}(t,s^{j},s^{\gra})\partial_{s^{\grb}}+\grh^{(1)i}_{\bar{l}}(t,s^{j},s^{\gra},\dot{s}^{j},\dot{s}^{\gra})\partial_{\dot{s}^{i}}+\grh^{(1)\grb}_{\bar{l}}(t,s^{j},s^{\gra},\dot{s}^{j},\dot{s}^{\gra})\partial_{\dot{s}^{\grb}}.\label{e28}
\end{align}
As can be deduced from \eqref{e27},\eqref{e28} we have not explicitly replace neither $\dot{s}^{i}$ by $u^{i}$, nor $s^{i}$ by $\int{u^{i}dt}+c^{i}$, the reason for this will become clear later. Furthermore, we have proven the condition (2) of the theorem as can be recognized by \eqref{e27}.

Given that prior to the reduction the following equations hold
\begin{align}
\left[\mathbf{Z}^{(1)}_{\bar{l}},\mathbf{A}\right]=-(\mathbf{A}(\grj_{\bar{l}}))\mathbf{A}\equiv\grl_{\bar{l}}\mathbf{A},\label{e29}
\end{align}
under which conditions the following equations are satisfied?
\begin{align}
\left[\mathbf{Y}_{\bar{l}},\hat{\mathbf{A}}\right]=-(\hat{\mathbf{A}}(\grj_{\bar{l}}))\hat{\mathbf{A}}\equiv\grm_{\bar{l}}\hat{\mathbf{A}}.\label{e30}
\end{align}
What we actually want to find is a relation between $Z_{i}$ and $Z_{\bar{l}}$ in order to prove the first point of the Theorem 1. Firstly we note that,
\begin{align}
&\mathbf{Z}^{(1)}_{\bar{l}}=\mathbf{Y}_{\bar{l}}+\grh^{l}_{\bar{l}}\mathbf{Z}^{(1)}_{l},\label{e313}\\
&\mathbf{A}=\hat{\mathbf{A}}+\dot{s}^{l}\mathbf{Z}^{(1)}_{l}.\label{e32}
\end{align}
If the above two expressions are used in \eqref{e29} we have that
\begin{align}
&\left[\mathbf{Y}_{\bar{l}}+\grh^{l}_{\bar{l}}\mathbf{Z}^{(1)}_{l},\hat{\mathbf{A}}+\dot{s}^{i}\mathbf{Z}^{(1)}_{i}\right]=\grl_{\bar{l}}(\hat{\mathbf{A}}+\dot{s}^{l}\mathbf{Z}^{(1)}_{l})\Rightarrow\nonumber\\
&\left[\mathbf{Y}_{\bar{l}},\hat{\mathbf{A}}\right]+\left[\grh^{l}_{\bar{l}}\mathbf{Z}^{(1)}_{l},\hat{\mathbf{A}}\right]+\left[\mathbf{Y}_{\bar{l}}+\grh^{l}_{\bar{l}}\mathbf{Z}^{(1)}_{l},\dot{s}^{l}\mathbf{Z}^{(1)}_{l}\right]=\grl_{\bar{l}}(\hat{\mathbf{A}}+\dot{s}^{i}\mathbf{Z}^{(1)}_{i})\Rightarrow\nonumber\\
&\left[\mathbf{Y}_{\bar{l}},\hat{\mathbf{A}}\right]+\left[\grh^{l}_{\bar{l}}\mathbf{Z}^{(1)}_{l},\hat{\mathbf{A}}\right]+\left[\mathbf{Y}_{\bar{l}}\left(\dot{s}^{i}\right)+\grh^{l}_{\bar{l}}\underbrace{\mathbf{Z}^{(1)}_{l}\left(\dot{s}^{i}\right)}_{=0}\right]\mathbf{Z}^{(1)}_{i}+\dot{s}^{i}\left[\mathbf{Y}_{\bar{l}}+\grh^{l}_{\bar{l}}\mathbf{Z}^{(1)}_{l},\mathbf{Z}^{(1)}_{i}\right]=\grl_{\bar{l}}(\hat{\mathbf{A}}+\dot{s}^{l}\mathbf{Z}^{(1)}_{l})\Rightarrow\nonumber\\
&\left[\mathbf{Y}_{\bar{l}},\hat{\mathbf{A}}\right]+\left[\grh^{l}_{\bar{l}}\mathbf{Z}^{(1)}_{l},\hat{\mathbf{A}}\right]+\grh^{(1)i}_{\bar{l}}\mathbf{Z}^{(1)}_{i}+\dot{s}^{i}\left[\mathbf{Z}^{(1)}_{\bar{l}},\mathbf{Z}^{(1)}_{i}\right]=\grl_{\bar{l}}(\hat{\mathbf{A}}+\dot{s}^{l}\mathbf{Z}^{(1)}_{l})\Rightarrow\nonumber\\
&\left[\mathbf{Y}_{\bar{l}},\hat{\mathbf{A}}\right]+\grh^{l}_{\bar{l}}\underbrace{\left[\mathbf{Z}^{(1)}_{l},\hat{\mathbf{A}}\right]}_{=0}-\hat{\mathbf{A}}\left(\grh^{l}_{\bar{l}}\right)\mathbf{Z}^{(1)}_{l}+\grh^{(1)i}_{\bar{l}}\mathbf{Z}^{(1)}_{i}-\dot{s}^{i}\left[\mathbf{Z}^{(1)}_{i},\mathbf{Z}^{(1)}_{\bar{l}}\right]=\grl_{\bar{l}}(\hat{\mathbf{A}}+\dot{s}^{l}\mathbf{Z}^{(1)}_{l})\Rightarrow\nonumber
\end{align}
\begin{align}
&\left[\mathbf{Y}_{\bar{l}},\hat{\mathbf{A}}\right]-\hat{\mathbf{A}}\left(\grh^{l}_{\bar{l}}\right)\mathbf{Z}^{(1)}_{l}+\grh^{(1)i}_{\bar{l}}\mathbf{Z}^{(1)}_{i}-\dot{s}^{i}\left[\mathbf{Z}^{(1)}_{i},\mathbf{Z}^{(1)}_{\bar{l}}\right]=\grl_{\bar{l}}(\hat{\mathbf{A}}+\dot{s}^{l}\mathbf{Z}^{(1)}_{l}).\label{e33}
\end{align}
At this point, if we want \eqref{e30} to hold, keeping in mind that
\begin{align}
\left[\mathbf{Z}^{(1)}_{i},\mathbf{Z}^{(1)}_{\bar{l}}\right]=C_{i\bar{l}}^{l}\mathbf{Z}^{(1)}_{l}+C^{\bar{k}}_{i\bar{l}}\mathbf{Z}^{(1)}_{\bar{k}}+C^{\grt}_{i\bar{l}}\mathbf{Z}^{(1)}_{\grt},\label{e34}
\end{align}
where $\mathbf{Z}^{(1)}_{\grt}$ the generators that are inherited, $\grt=\bar{m}+1,..,r$, then equation \eqref{e33} becomes
\begin{align}
\left(\grm_{\bar{l}}-\grl_{\bar{l}}\right)\hat{\mathbf{A}}+\left[-\hat{\mathbf{A}}\left(\grh^{l}_{\bar{l}}\right)+\grh^{(1)l}_{\bar{l}}-\dot{s}^{i}C^{l}_{i\bar{l}}-\grl_{\bar{l}}\dot{s}^{l}\right]\mathbf{Z}^{(1)}_{l}-\dot{s}^{i}C^{\bar{k}}_{i\bar{l}}\mathbf{Z}^{(1)}_{\bar{k}}-\dot{s}^{i}C^{\grt}_{i\bar{l}}\mathbf{Z}^{(1)}_{\grt}=0.\label{e35}
\end{align}
By assumption, $\hat{A}, \mathbf{Z}^{(1)}_{l},\mathbf{Z}^{(1)}_{\bar{l}},\mathbf{Z}^{(1)}_{\grt}$ are linearly independent, thus the equation \eqref{e35} splits into four sets of equations
\begin{align}
&\grm_{\bar{l}}=\grl_{\bar{l}},\label{e36}\\
&-\hat{\mathbf{A}}\left(\grh^{l}_{\bar{l}}\right)+\grh^{(1)l}_{\bar{l}}-\dot{s}^{i}C^{l}_{i\bar{l}}=\grl_{\bar{l}}\dot{s}^{l},\label{e37}\\
&\dot{s}^{i}C^{\bar{k}}_{i\bar{l}}=0,\label{e38}\\
&\dot{s}^{i}C^{\grt}_{i\bar{l}}=0.\label{e38887}
\end{align}   
The third and fourth set $\eqref{e38}$, $\eqref{e38887}$ due to the linear independence of $\dot{s}^{i}$ imply
\begin{align}
C^{\bar{k}}_{i\bar{l}}=0,\label{e39}\\
C^{\grt}_{i\bar{l}}=0,\label{e39trtr}
\end{align}
which in combination with \eqref{e34} leads to the desired condition
\begin{align}
&\left[\mathbf{Z}^{(1)}_{i},\mathbf{Z}^{(1)}_{\bar{l}}\right]=C_{i\bar{l}}^{l}\mathbf{Z}^{(1)}_{l}\Leftrightarrow\nonumber\\
&\left[\mathbf{Z}_{i},\mathbf{Z}_{\bar{l}}\right]=C_{i\bar{l}}^{l}\mathbf{Z}_{l},\label{e40}
\end{align}
as it is expressed in condition (1) of the theorem.

The other two sets give us information about the form of the remaining generators \eqref{e28}.
\begin{align}
&\grm_{\bar{l}}=\grl_{\bar{l}}\Rightarrow\nonumber\\
&-\hat{\mathbf{A}}(\grj_{\bar{l}})=-\mathbf{A}(\grj_{\bar{l}})\Rightarrow\nonumber\\
&\left(\mathbf{A}-\hat{\mathbf{A}}\right)\left(\grj_{\bar{l}}\right)=0\Rightarrow\nonumber\\
&\dot{s}^{l}\mathbf{Z}^{(1)}_{l}\left(\grj_{\bar{l}}\right)=0\Rightarrow\nonumber\\
&\mathbf{Z}^{(1)}_{l}\left(\grj_{\bar{l}}\right)=0\Rightarrow\nonumber\\
&\partial_{s^{l}}\grj_{\bar{l}}(t,s^{j},s^{\gra})=0\Rightarrow\nonumber
\end{align}
\begin{align}
\grj_{\bar{l}}(t,s^{j},s^{\gra})=\grj_{\bar{l}}(t,s^{\gra}).\label{e41}
\end{align}
From \eqref{e37}, by recalling an alternative, equivalent definition of $\grh^{(1)l}_{\bar{l}}(t,s^{j},s^{\gra},\dot{s}^{j},\dot{s}^{\gra})$ namely
\begin{align}
\grh^{(1)l}_{\bar{l}}(t,s^{j},s^{\gra},\dot{s}^{j},\dot{s}^{\gra})=\mathbf{A}\left[\grh_{\bar{l}}^{l}(t,s^{j},s^{\gra})\right]+\dot{s}^{l}\grl_{\bar{l}},\label{e42}
\end{align}
we get
\begin{align}
&-\hat{\mathbf{A}}\left(\grh^{l}_{\bar{l}}\right)+\mathbf{A}\left(\grh_{\bar{l}}^{l}\right)+\dot{s}^{l}\grl_{\bar{l}}-\dot{s}^{i}C^{l}_{i\bar{l}}=\grl_{\bar{l}}\dot{s}^{l}\Rightarrow\nonumber\\
&\left(\mathbf{A}-\hat{\mathbf{A}}\right)\left(\grh^{l}_{\bar{l}}\right)-\dot{s}^{i}C_{i\bar{l}}^{l}=0\Rightarrow\nonumber\\
&\dot{s}^{i}\mathbf{Z}^{(1)}_{i}\left(\grh^{l}_{\bar{l}}\right)=\dot{s}^{i}C_{i\bar{l}}^{l}\Rightarrow\nonumber\\
&\mathbf{Z}^{(1)}_{i}\left(\grh^{l}_{\bar{l}}\right)=C_{i\bar{l}}^{l}\Rightarrow\nonumber\\
&\partial_{s^{i}}\grh^{l}_{\bar{l}}(t,s^{j},s^{\gra})=C_{i\bar{l}}^{l}\Rightarrow\nonumber
\end{align}
\begin{align}
\grh^{l}_{\bar{l}}(t,s^{j},s^{\gra})=s^{j}C_{j\bar{l}}^{l}+\bar{\grh}^{l}_{\bar{l}}(t,s^{\gra}).\label{e43}
\end{align}
Next we use the explicit forms of the prolongations in the equation \eqref{e40} in order to find some conditions on the set of functions $\grh^{\grb}_{\bar{l}}(t,s^{j},s^{\gra}),\grh^{(1)i}_{\bar{l}}(t,s^{j},s^{\gra},\dot{s}^{j},\dot{s}^{\gra}), \grh^{(1)\grb}_{\bar{l}}(t,s^{j},s^{\gra},\dot{s}^{j},\dot{s}^{\gra})$.
\begin{align}
&\left[\mathbf{Z}^{(1)}_{i},\mathbf{Z}^{(1)}_{\bar{l}}\right]=C_{i\bar{l}}^{l}\mathbf{Z}^{(1)}_{l}\Rightarrow\nonumber\\
&\underbrace{\partial_{s^{i}}\grj_{\bar{l}}}_{=0}\partial_{t}+\underbrace{\partial_{s^{i}}\grh^{l}_{\bar{l}}}_{=C^{l}_{i\bar{l}}}\partial_{s^{l}}+\partial_{s^{i}}\grh^{\grb}_{\bar{l}}\partial_{s^{\grb}}+\partial_{s^{i}}\grh^{(1)l}_{\bar{l}}\partial_{\dot{s}^{l}}+\partial_{s^{i}}\grh^{(1)\grb}_{\bar{l}}\partial_{\dot{s}^{\grb}}=C^{l}_{i\bar{l}}\partial_{s^{l}}\Rightarrow\nonumber\\
&C^{l}_{i\bar{l}}\partial_{s^{l}}+\partial_{s^{i}}\grh^{\grb}_{\bar{l}}\partial_{s^{\grb}}+\partial_{s^{i}}\grh^{(1)l}_{\bar{l}}\partial_{\dot{s}^{l}}+\partial_{s^{i}}\grh^{(1)\grb}_{\bar{l}}\partial_{\dot{s}^{\grb}}=C^{l}_{i\bar{l}}\partial_{s^{l}}\Rightarrow\nonumber\\
&\partial_{s^{i}}\grh^{\grb}_{\bar{l}}\partial_{s^{\grb}}+\partial_{s^{i}}\grh^{(1)l}_{\bar{l}}\partial_{\dot{s}^{l}}+\partial_{s^{i}}\grh^{(1)\grb}_{\bar{l}}\partial_{\dot{s}^{\grb}}=0\Rightarrow\nonumber
\end{align}
due to the independence of the basis $\partial_{s^{\grb}},\partial_{\dot{s}^{l}},\partial_{\dot{s}^{\grb}}$
\begin{align}
&\partial_{s^{i}}\grh^{\grb}_{\bar{l}}(t,s^{j},s^{\gra})=0,\label{e44}\\
&\partial_{s^{i}}\grh^{(1)l}_{\bar{l}}(t,s^{j},s^{\gra},\dot{s}^{j},\dot{s}^{\gra})=0,\label{e45}\\
&\partial_{s^{i}}\grh^{(1)\grb}_{\bar{l}}(t,s^{j},s^{\gra},\dot{s}^{j},\dot{s}^{\gra})=0.\label{e46}
\end{align}
Thus, the final form of the reduced inherited generators read
\begin{align}
&\mathbf{Y}_{\bar{l}}=\grj_{\bar{l}}(t,s^{\gra})\partial_{t}+\grh^{\grb}_{\bar{l}}(t,s^{\gra})\partial_{s^{\grb}}+\grh^{(1)i}_{\bar{l}}(t,s^{\gra},u^{j},\dot{s}^{\gra})\partial_{u^{i}}+\grh^{(1)\grb}_{\bar{l}}(t,s^{\gra},u^{j},\dot{s}^{\gra})\partial_{\dot{s}^{\grb}},\label{e48}
\end{align}
and thus, this proves the condition (3). For the fourth and final condition we start with the expression for the algebra of the generators to be inherited. In general
\begin{align}
[\mathbf{Z}^{(1)}_{\bar{l}},\mathbf{Z}^{(1)}_{\bar{k}}]=C^{\bar{q}}_{\bar{l}\bar{k}}\mathbf{Z}^{(1)}_{\bar{q}}+C^{l}_{\bar{l}\bar{k}}\mathbf{Z}^{(1)}_{l}+C^{\grt}_{\bar{l}\bar{k}}\mathbf{Z}^{(1)}_{\grt}.
\end{align}
Let us recall the relation \eqref{e313} and use it in the above equation
\begin{align}
&[\mathbf{Y}_{\bar{l}}+\grh^{l}_{\bar{l}}\mathbf{Z}_{l}^{(1)},\mathbf{Y}_{\bar{k}}+\grh^{j}_{\bar{k}}\mathbf{Z}_{j}^{(1)}]=C^{\bar{q}}_{\bar{l}\bar{k}}\mathbf{Y}_{\bar{q}}+C^{\bar{q}}_{\bar{l}\bar{k}}\grh_{\bar{q}}^{l}\mathbf{Z}^{(1)}_{l}+C^{l}_{\bar{l}\bar{k}}\mathbf{Z}^{(1)}_{l}+C^{\grt}_{\bar{l}\bar{k}}\mathbf{Z}^{(1)}_{\grt}\Rightarrow\nonumber\\
&[\mathbf{Y}_{\bar{l}},\mathbf{Y}_{\bar{k}}]+\mathbf{Y}_{\bar{l}}\left(\grh_{\bar{k}}^{l}\right)\mathbf{Z}^{(1)}_{l}+\grh^{l}_{\bar{k}}\underbrace{[\mathbf{Y}_{\bar{l}},\mathbf{Z}^{(1)}_{l}]}_{=0}+\grh^{l}_{\bar{l}}\underbrace{[\mathbf{Z}^{(1)}_{l},\mathbf{Y}_{\bar{k}}]}_{=0}-\mathbf{Y}_{\bar{k}}\left(\grh^{l}_{\bar{l}}\right)\mathbf{Z}^{(1)}_{l}+\grh_{\bar{l}}^{j}\mathbf{Z}^{(1)}_{j}\left(\grh_{\bar{k}}^{l}\right)\mathbf{Z}^{(1)}_{l}-\grh_{\bar{k}}^{j}\mathbf{Z}^{(1)}_{j}\left(\grh_{\bar{l}}^{l}\right)\mathbf{Z}^{(1)}_{l}\nonumber\\
&+\grh^{l}_{\bar{l}}\grh^{j}_{\bar{k}}\underbrace{[\mathbf{Z}^{(1)}_{l},\mathbf{Z}^{(1)}_{j}]}_{=0}=C^{\bar{q}}_{\bar{l}\bar{k}}\mathbf{Y}_{\bar{q}}+C^{\bar{q}}_{\bar{l}\bar{k}}\grh_{\bar{q}}^{l}\mathbf{Z}^{(1)}_{l}+C^{l}_{\bar{l}\bar{k}}\mathbf{Z}^{(1)}_{l}+C^{\grt}_{\bar{l}\bar{k}}\mathbf{Z}^{(1)}_{\grt}\Rightarrow\nonumber\\
&[\mathbf{Y}_{\bar{l}},\mathbf{Y}_{\bar{k}}]-C^{\bar{q}}_{\bar{l}\bar{k}}\mathbf{Y}_{\bar{q}}+\left[\mathbf{Y}_{\bar{l}}\left(\grh_{\bar{k}}^{l}\right)-\mathbf{Y}_{\bar{k}}\left(\grh^{l}_{\bar{l}}\right)+\grh_{\bar{l}}^{j}\mathbf{Z}^{(1)}_{j}\left(\grh_{\bar{k}}^{l}\right)-\grh_{\bar{k}}^{j}\mathbf{Z}^{(1)}_{j}\left(\grh_{\bar{l}}^{l}\right)-C^{\bar{q}}_{\bar{l}\bar{k}}\grh_{\bar{q}}^{l}-C^{l}_{\bar{l}\bar{k}}\right]\mathbf{Z}^{(1)}_{l}-C^{\grt}_{\bar{l}\bar{k}}\mathbf{Z}^{(1)}_{\grt}=0.
\end{align}
Given that $\mathbf{Z}^{(1)}_{l}, \mathbf{Z}^{(1)}_{\grt}$ are independent themselves and independent from $\mathbf{Y}_{\bar{q}}, [\mathbf{Y}_{\bar{l}},\mathbf{Y}_{\bar{k}}]$ we are led to the equations
\begin{align}
&[\mathbf{Y}_{\bar{l}},\mathbf{Y}_{\bar{k}}]=C^{\bar{q}}_{\bar{l}\bar{k}}\mathbf{Y}_{\bar{q}},\label{545442}\\
&\mathbf{Y}_{\bar{l}}\left(\grh_{\bar{k}}^{l}\right)-\mathbf{Y}_{\bar{k}}\left(\grh^{l}_{\bar{l}}\right)+\grh_{\bar{l}}^{j}\mathbf{Z}^{(1)}_{j}\left(\grh_{\bar{k}}^{l}\right)-\grh_{\bar{k}}^{j}\mathbf{Z}^{(1)}_{j}\left(\grh_{\bar{l}}^{l}\right)-C^{\bar{q}}_{\bar{l}\bar{k}}\grh_{\bar{q}}^{l}-C^{l}_{\bar{l}\bar{k}}=0,\label{sdsdssa}\\
&C^{\grt}_{\bar{l}\bar{k}}=0.\label{aacddf}
\end{align}
The equation \eqref{545442} proves the condition (4). The equation \eqref{aacddf} give us additional information about the structure constants related to the generators that will not be inherited. This additional piece of information lead us to the following corollary, by taking into account \eqref{e11}, \eqref{e39trtr}, \eqref{aacddf}:

\textbf{\underline{Corollary 1}:}
\textit{\textbf{The subset $\{\mathbf{Z}_{l},\mathbf{Z}_{\bar{l}}\}$ of the set $\{\mathbf{Z}_{l},\mathbf{Z}_{\bar{l}}, \mathbf{Z}_{\grt}\}$ forms a sub-algebra of the original algebra:}}
\begin{align}
&[\mathbf{Z}_{i},\mathbf{Z}_{j}]=0,\label{etrg3}\\
&[\mathbf{Z}_{i},\mathbf{Z}_{\bar{l}}]=C^{l}_{i\bar{l}}\mathbf{Z}_{l},\label{etrdgh1}\\
&[\mathbf{Z}_{\bar{l}},\mathbf{Z}_{\bar{k}}]=C^{l}_{\bar{l}\bar{k}}\mathbf{Z}_{l}+C^{\bar{q}}_{\bar{l}\bar{k}}\mathbf{Z}_{\bar{q}}.\label{egggh2}
\end{align} 
When it comes to the equation \eqref{sdsdssa}, it restricts even more the form for the components of $\mathbf{Y}_{\bar{l}}$.

\section{Solvable algebra}

For the purposes of this section, we follow the notation of \cite{Gilmore:102082}. Let us start by giving the definition of an invariant sub-algebra.

\textbf{\underline{Definition 1}:} \textit{\textbf{An invariant sub-algebra h of an algebra g is an algebra whose the commutator of it's every element with every element of g, belongs to the sub-algebra. Symbolically this is written as}}
\begin{align}
\left[\mathbf{h},\mathbf{g}\right]\subseteq \mathbf{h}.\label{e58}
\end{align}  

\textbf{\underline{Definition 2}:} \textit{\textbf{An algebra will be called (n)-level solvable if it admits a series of invariant sub-algebras such that}}
\begin{align}
\mathbf{g}\equiv\mathbf{g}^{(0)}\supset\mathbf{g}^{(1)}\supset\mathbf{g}^{(2)}\supset...\supset\mathbf{g}^{(n-1)}\supset\mathbf{g}^{(n)}\equiv\left\{0\right\}.\label{e59}
\end{align}
A way to always construct invariant sub-algebras is by considering the derived algebra.

\textbf{\underline{Definition 3}:} \textit{\textbf{Derived algebra is called an algebra $\mathbf{g}^{(1)}$ which is constructed by some linearly independent sub-set of elements of the commutator of the algebra $\mathbf{g}^{(0)}\equiv \mathbf{g}$. Symbolically,}}
\begin{align}
\left[\mathbf{g}^{(0)},\mathbf{g}^{(0)}\right]\equiv \mathbf{g}^{(1)},\label{e60}
\end{align}
\textbf{\textit{Ascribable to this definition, the following conditions hold}}
\begin{align}
&\mathbf{g}^{(1)}\subseteq \mathbf{g}^{(0)},\label{e61}\\
&\left[\mathbf{g}^{(1)},\mathbf{g}^{(0)}\right]\subseteq\mathbf{g}^{(1)}.\label{e62}
\end{align}
Thus, due to \eqref{e62}, the sub-algebra $\mathbf{g}^{(1)}$ is an invariant sub-algebra of $\mathbf{g}^{(0)}$. Let us now assume a (3)-level solvable algebra which is constructed by the derived algebras
\begin{align}
\mathbf{g}\equiv\mathbf{g}^{(0)}\supset\mathbf{g}^{(1)}\supset\mathbf{g}^{(2)}\supset\mathbf{g}^{(3)}\equiv\left\{0\right\}.\label{e63}
\end{align}
The last non-trivial sub-algebra is $\mathbf{g}^{(2)}$ which also has the property of being Abelian. Due to the definition of the derived algebras the following relations hold
\begin{align}
\left[\mathbf{g}^{(0)},\mathbf{g}^{(0)}\right]\equiv\mathbf{g}^{(1)},\hspace{0.2cm}\left[\mathbf{g}^{(1)},\mathbf{g}^{(1)}\right]\equiv\mathbf{g}^{(2)},\hspace{0.2cm}\left[\mathbf{g}^{(1)},\mathbf{g}^{(0)}\right]\subseteq\mathbf{g}^{(1)},\hspace{0.2cm}\left[\mathbf{g}^{(2)},\mathbf{g}^{(1)}\right]\subseteq\mathbf{g}^{(2)}.\hspace{0.2cm}\label{e64}
\end{align}
As can be observed from \eqref{e64}, the algebra $\mathbf{g}^{(1)}$ is an invariant sub-algebra of $\mathbf{g}^{(0)}$ and $\mathbf{g}^{(2)}$ is an invariant sub-algebra of $\mathbf{g}^{(1)}$. It is not difficult to prove that due to the Jacobi identity $\mathbf{g}^{(2)}$ is also an invariant sub-algebra of $\mathbf{g}^{(0)}$. Symbolically we can write this as follows
\begin{align}
&\left[\mathbf{g}^{(2)},\mathbf{g}^{(0)}\right]\equiv\left[\left[\mathbf{g}^{(1)},\mathbf{g}^{(1)}\right],\mathbf{g}^{(0)}\right]=-\left[\mathbf{g}^{(0)},\left[\mathbf{g}^{(1)},\mathbf{g}^{(1)}\right]\right]\xRightarrow{\text{Jacobi identity}}\nonumber\\
&\left[\mathbf{g}^{(2)},\mathbf{g}^{(0)}\right]=\left[\mathbf{g}^{(1)},\underbrace{\left[\mathbf{g}^{(1)},\mathbf{g}^{(0)}\right]}_{\subseteq\mathbf{g}^{(1)}}\right]+\left[\mathbf{g}^{(1)},\underbrace{\left[\mathbf{g}^{(0)},\mathbf{g}^{(1)}\right]}_{\subseteq\mathbf{g}^{(1)}}\right]\Rightarrow\nonumber\\
&\left[\mathbf{g}^{(2)},\mathbf{g}^{(0)}\right]\subseteq\left[\mathbf{g}^{(1)},\mathbf{g}^{(1)}\right]\equiv\mathbf{g}^{(2)}\Rightarrow\nonumber\\
&\left[\mathbf{g}^{(2)},\mathbf{g}^{(0)}\right]\subseteq\mathbf{g}^{(2)}.\label{e65}
\end{align}
This can be generalized to any (n)-level solvable algebra, thus we arrive at the following Theorem.

\textbf{\underline{Theorem 2}:} \textbf{\textit{For all (n)-level solvable algebras the following relations hold}}
\begin{align}
\left[\mathbf{g}^{(i)},\mathbf{g}^{(j)}\right]\subseteq\mathbf{g}^{(k)},\hspace{0.2cm} k=max(i,j),\hspace{0.2cm} i\neq{j},\label{e66}
\end{align}
\textbf{\textit{where $\mathbf{g}^{(l)}$ is the derived algebra $\mathbf{g}^{(l)}=\left[\mathbf{g}^{(l-1)},\mathbf{g}^{(l-1)}\right]$ and $max(i,j)$ the maximal between the indices $i,j$.}} 

From this theorem, the following corollary can be deduced:

\textbf{\underline{Corollary 2}:} \textbf{\textit{Every (n)-level solvable algebra consists of (n)-number invariant sub-algebras, with the (n)-th being the empty set.}}
 
At last, let us define a ``coset". 

\textbf{\underline{Definition 4}:} \textbf{\textit{The difference between two derived algebras $\mathbf{g}^{(i)}$, $\mathbf{g}^{(j)}$ will be called the ``coset"}}
\begin{align}
\mathbf{B}^{(i)}_{(j)}=\mathbf{g}^{(i)}-\mathbf{g}^{(j)},\label{e67}
\end{align} 
\textbf{where $i<j$.}

\section{Combination of solvability and ``first integration method''}

In this section we would like to combine the knowledge of the two previous sections. To prove the theorem below, we have to interplay between the abstract notation for $\mathbf{g}^{(i)}$, $\mathbf{B}^{(i)}_{(j)}$ and the structure constants.

\textbf{\underline{Theorem 3}:} \textit{\textbf{For an (n)-level solvable algebra $\mathbf{g}^{(0)}$, by applying the ``first integration method'' on the invariant derived sub-algebra $\mathbf{g}^{(n-1)}$, it is guaranteed that the prolonged and reduced generators of the ``coset" $\mathbf{B}^{(0)}_{(n-1)}$ will be inherited on the reduced equations.}}

\textbf{\underline{Proof}}

We can easily prove the above theorem based on previous theorems and definitions. We start from the fact that when there is an (n)-level solvable algebra $\mathbf{g}^{(0)}$, according to the definitions of solvable and derived algebras, there is certainly an Abelian algebra which is provided by the derived algebra $\mathbf{g}^{(n-1)}$. Furthermore, due to Theorem 2, this is also an invariant sub-algebra, so based on \eqref{e66}, $\left[\mathbf{g}^{(n-1)},\mathbf{g}^{(0)}\right]\subseteq\mathbf{g}^{(n-1)}$. From this relation we infer that
\begin{align}
\left[\mathbf{g}^{(n-1)},\mathbf{B}^{(0)}_{(n-1)}\right]\subseteq\mathbf{g}^{(n-1)},\label{aq1}
\end{align}
where $\mathbf{B}^{(0)}_{(n-1)}=\mathbf{g}^{(0)}-\mathbf{g}^{(n-1)}$.
The equation \eqref{aq1} is the abstract notation of the equation $\mathbf{(1)}$ which is stated in Theorem 1 and provide us with the relation that a generator needs to satisfy in order for it's first prolongation to be inherited as a symmetry of the reduced system, after the use of the first integration method on the Abelian sub-algebra $\mathbf{g}^{(n-1)}$. This concludes the proof.

So far, solvability provide us with the ``best'' Abelian sub-algebra on which we have to apply the first integration method, in order to inherit the maximal number of the generators in the reduced equations. The question that we like to address now is the following:

Do the prolonged and reduced elements of the ``coset'' $\mathbf{B}^{(0)}_{(n-1)}$ form again a solvable algebra so that we could repeat the procedure? 

The answer to this question is providing by the following theorem

\textbf{\underline{Theorem 4}:} \textit{\textbf{Given an (n)-level solvable algebra, the first integration procedure can be repeated (n)-times upon the following chain of ``cosets" }}
\begin{align}
&g^{(n-1)}\equiv\mathbf{B}^{(n-1)}_{(n)}\xrightarrow{pr}\mathbf{B}^{(n-2)}_{(n-1)}\xrightarrow{pr}....\xrightarrow{pr}\mathbf{B}^{(0)}_{(1)},\label{aq2}
\end{align}
\textit{\textbf{starting from $\mathbf{B}^{(n-1)}_{(n)}$, given that the generators of each ``coset'' act transitively in some proper subspace. The symbol $\mathbf{pr}$ over the arrow implies that at each ``coset" except the $\mathbf{B}^{(n-1)}_{(n)}$, we have to consider the proper prolonged and reduced form of it's generators.}}

We will prove this theorem in an iterative sense. Before we do this, it is instructive to provide the following very important Corollary of the previous theorem.  

\textbf{\underline{Corollary 3}:} \textbf{\textit{For an (n)-level solvable algebra, all the generators can be used via the ``first integration method''.}}

Let us start with a (2)-level solvable algebra. In some cases we are going to use tables, in order for the various notations to be more transparent.

\subsection{(2)-level solvable algebra}

Due to the definition of solvablility we get
\begin{align}
\left\{0\right\}\equiv\mathbf{g}^{(2)}\subset\mathbf{g}^{(1)}\subset\mathbf{g}^{(0)}.\label{be1}
\end{align}
Furthermore, due to \eqref{e66}
\begin{align}
&\left[\mathbf{g}^{(0)},\mathbf{g}^{(0)}\right]\equiv\mathbf{g}^{(1)}\Rightarrow \left[\mathbf{g}^{(1)},\mathbf{g}^{(0)}\right]\subseteq\mathbf{g}^{(1)}, \left[\mathbf{B}^{(0)}_{(1)},\mathbf{g}^{(0)}\right]\subseteq\mathbf{g}^{(1)}, \left[\mathbf{B}^{(0)}_{(1)},\mathbf{B}^{(0)}_{(1)}\right]\subseteq\mathbf{g}^{(1)},\label{be2}\\
&\left[\mathbf{g}^{(1)},\mathbf{g}^{(1)}\right]\equiv\mathbf{g}^{(2)}\equiv\left\{0\right\}\Rightarrow \left[\left\{0\right\},\mathbf{g}^{(1)}\right]\subseteq\left\{0\right\}, \left[\mathbf{B}^{(1)}_{(2)},\mathbf{g}^{(0)}\right]\subseteq\left\{0\right\}.\label{be3}
\end{align}
\subsubsection{Indices}
Because we are going to translate the previous relations into structure constants, we have to define some indices.
\begin{center}
\begin{tabular}{ |c|c|c| } 
\hline
 \textbf{Indices of $\mathbf{g}^{(i)}$} & \textbf{Indices of $\mathbf{B}^{(i)}_{(j)}$}  & \textbf{Useful relations} \\
\hline
$\grm_{0},\grn_{0},\grl_{0},..\rightarrow \mathbf{g}^{(0)}$ & $i_{1},j_{1},l_{1},..\rightarrow \mathbf{B}^{(0)}_{(1)}$ & $\grm_{0}=\grm_{1}+i_{1}$\\
\hline
$\grm_{1},\grn_{1},\grl_{1},..\rightarrow \mathbf{g}^{(1)}$ & $i_{2},j_{2},l_{2},..\rightarrow \mathbf{B}^{(1)}_{(2)}$ & $\grm_{1}=\grm_{2}+i_{2}$ \\
\hline
$\grm_{2},\grn_{2},\grl_{2},..\rightarrow \mathbf{g}^{(2)}$ & $i_{3},j_{3},l_{3},..\rightarrow \mathbf{B}^{(0)}_{(2)}$ & $i_{3}=i_{1}+i_{2}=\grm_{0}$ \\
\hline
\end{tabular}
\captionof{table}{The expressions of the third column can be deduced from the definition of the ``coset" and the derived algebra. Also, they hold for all the indices of the same family.}
\end{center} 
\subsubsection{The structure constants}
The structure constants of the algebra $g^{(0)}$ will be symbolized as $C^{\grl_{0}}_{\grm_{0}\grn_{0}}$. From the \eqref{be2},\eqref{be3} we get
\begin{align}
&\eqref{be2}\Rightarrow C^{l_{1}}_{\grm_{0}\grn_{0}}=0,\label{be4}\\
&\eqref{be3}\Rightarrow C^{l_{3}}_{\grm_{1}\grn_{1}}=0\equiv C^{\grl_{0}}_{\grm_{1}\grn_{1}}=0.\label{be5}
\end{align}
\subsubsection{The reduced system}
Due to the fact that $\mathbf{B}^{(1)}_{(2)}=g^{(1)}-g^{(2)}=g^{(1)}$, $\mathbf{B}^{(1)}_{(2)}$ forms an Abelian algebra, we can use the first integration method upon it. The symmetries to be inherited will be given by the prolongation of the generators of $B^{(0)}_{(1)}$ properly reduced. As we have learned from $(\mathbf{4})$ of Theorem 1, the reduced generators will form an algebra with structure constants a subset of the original ones, which in our notation read $C^{l_{1}}_{i_{1}j_{1}}$. Of course, since the structure constants does not change, all the previous equations \eqref{be4},\eqref{be5} could still provide us with useful information. We infer, by using the expressions from the third column of Table 1, that
\begin{align}
\eqref{be4}\Rightarrow C^{l_{1}}_{i_{1}j_{1}}=0\Leftrightarrow \left[\mathbf{B}^{(0)}_{(1)},\mathbf{B}^{(0)}_{(1)}\right]\equiv\left\{0\right\}.\label{be6}
\end{align}  
Note that this happens only for the properly reduced generators of $\mathbf{B}^{(0)}_{(1)}$. As a result, we can apply the ``first integration method'' once more on the elements of the algebra $\mathbf{B}^{(0)}_{(1)}$ and arriving at the following chain
\begin{align}
\mathbf{B}^{(1)}_{(2)}\xrightarrow{pr}\mathbf{B}^{(0)}_{(1)}.\label{be7}
\end{align}

\subsection{(3)-level solvable algebra}

This case is more evolved but nevertheless, the procedure is the same.
\begin{align}
\left\{0\right\}\equiv\mathbf{g}^{(3)}\subset\mathbf{g}^{(2)}\subset\mathbf{g}^{(1)}\subset\mathbf{g}^{(0)},\label{be8}
\end{align}
while from \eqref{e66}
\begin{align}
&\left[\mathbf{g}^{(0)},\mathbf{g}^{(0)}\right]\equiv \mathbf{g}^{(1)}\Rightarrow \left[\mathbf{g}^{(1)},\mathbf{g}^{(0)}\right]\subseteq\mathbf{g}^{(1)}, \left[\mathbf{B}^{(0)}_{(1)},\mathbf{g}^{(0)}\right]\subseteq\mathbf{g}^{(1)}, \left[\mathbf{B}^{(0)}_{(1)},\mathbf{B}^{(0)}_{(1)}\right]\subseteq\mathbf{g}^{(1)},\label{be9}\\
&\left[\mathbf{g}^{(1)},\mathbf{g}^{(1)}\right]\equiv \mathbf{g}^{(2)}\Rightarrow \left[\mathbf{g}^{(2)},\mathbf{g}^{(1)}\right]\subseteq\mathbf{g}^{(2)}, \left[\mathbf{B}^{(1)}_{(2)},\mathbf{g}^{(1)}\right]\subseteq\mathbf{g}^{(2)}, \left[\mathbf{B}^{(1)}_{(2)},\mathbf{B}^{(1)}_{(2)}\right]\subseteq\mathbf{g}^{(2)},\label{be10}\\
&\left[\mathbf{g}^{(2)},\mathbf{g}^{(2)}\right]\equiv \mathbf{g}^{(3)}\equiv\left\{0\right\}\Rightarrow \left[\left\{0\right\},\mathbf{g}^{(2)}\right]\subseteq\left\{0\right\}, \left[\mathbf{B}^{(2)}_{(3)},\mathbf{g}^{(2)}\right]\subseteq\left\{0\right\},\label{be11}\\
&\left[\mathbf{g}^{(2)},\mathbf{g}^{(0)}\right]\subseteq\mathbf{g}^{(2)}.\label{be12}
\end{align}
\subsubsection{Indices}
The corresponding table now reads
\begin{center}
\begin{tabular}{ |c|c|c| } 
\hline
 \textbf{Indices of $\mathbf{g}^{(i)}$} & \textbf{Indices of $\mathbf{B}^{(i)}_{(j)}$}  & \textbf{Useful relations} \\
\hline
$\grm_{0},\grn_{0},\grl_{0},..\rightarrow \mathbf{g}^{(0)}$ & $i_{1},j_{1},l_{1},..\rightarrow \mathbf{B}^{(0)}_{(1)}$ & $\grm_{0}=\grm_{1}+i_{1}$\\
\hline
$\grm_{1},\grn_{1},\grl_{1},..\rightarrow \mathbf{g}^{(1)}$ & $i_{2},j_{2},l_{2},..\rightarrow \mathbf{B}^{(1)}_{(2)}$ & $\grm_{1}=\grm_{2}+i_{2}$ \\
\hline
$\grm_{2},\grn_{2},\grl_{2},..\rightarrow \mathbf{g}^{(2)}$ & $i_{3},j_{3},l_{3},..\rightarrow \mathbf{B}^{(0)}_{(2)}$ & $i_{3}=i_{1}+i_{2}$ \\
\hline
$\grm_{3},\grn_{3},\grl_{3},..\rightarrow \mathbf{g}^{(3)}$ & $i_{4},j_{4},l_{4},..\rightarrow \mathbf{B}^{(2)}_{(3)}$ & $i_{4}=\grm_{2}$ \\
\hline
 & $i_{5},j_{5},l_{5},..\rightarrow \mathbf{B}^{(1)}_{(3)}$ & $i_{5}=i_{2}+i_{4}=\grm_{1}$ \\
\hline
 & $i_{6},j_{6},l_{6},..\rightarrow \mathbf{B}^{(0)}_{(3)}$ & $i_{6}=i_{1}+i_{2}+i_{4}=i_{3}+i_{4}=i_{3}+\grm_{2}=i_{1}+\grm_{1}=\grm_{0}$ \\
\hline
\end{tabular}
\captionof{table}{The expressions of the third column can be deduced from the definition of the ``coset" and the derived algebra}
\end{center} 
\subsubsection{The structure constants}
The structure constants of $\mathbf{g^{(0)}}$ will symbolized as $C^{\grl_{0}}_{\grm_{0}\grn_{0}}$.
\begin{align}
&\eqref{be9}\Rightarrow C^{l_{1}}_{\grm_{0}\grn_{0}}=0,\label{be13}\\
&\eqref{be10}\Rightarrow C^{l_{3}}_{\grm_{1}\grn_{1}}=0,\label{be14}\\
&\eqref{be11}\Rightarrow C^{l_{6}}_{\grm_{2}\grn_{2}}=0\equiv C^{\grl_{0}}_{\grm_{2}\grn_{2}}=0,\label{be15}\\
&\eqref{be12}\Rightarrow C^{l_{3}}_{\grm_{2}\grn_{0}}=0.\label{be16}
\end{align}
\subsubsection{The reduced system}
Same as before, $\mathbf{B}^{(2)}_{(3)}=g^{(2)}-g^{(3)}=g^{(2)}$ forms an Abelian invariant sub-algebra. The symmetries to be inherited will consist of the prolonged elements of the ``coset" $\mathbf{B}^{(0)}_{(2)}$. The corresponding structure constants will be $C^{l_{3}}_{i_{3}j_{3}}$. From the third row of the third column, we can split those structure constants into $C^{l_{1}}_{i_{3}j_{3}}$, $C^{l_{2}}_{i_{3}j_{3}}$, since $l_{3}=l_{1}+l_{2}$. Due to \eqref{be13} and the final row of the third column we can deduce that $C^{l_{1}}_{i_{3}j_{3}}=0$. This implies that,
\begin{align}
\left[\mathbf{B}^{(0)}_{(2)},\mathbf{B}^{(0)}_{(2)}\right]\equiv\mathbf{B}^{(1)}_{(2)}\Rightarrow\left[\mathbf{B}^{(1)}_{(2)},\mathbf{B}^{(0)}_{(2)}\right]\subseteq\mathbf{B}^{(1)}_{(2)},\left[\mathbf{B}^{(0)}_{(1)},\mathbf{B}^{(0)}_{(2)}\right]\subseteq\mathbf{B}^{(1)}_{(2)}.\label{be17}
\end{align}
Since $\mathbf{B}^{(1)}_{(2)}$ is the derived algebra, we may consider also it's derived algebra. In general, since is an invariant sub-algebra we expect that $\left[\mathbf{B}^{(1)}_{(2)},\mathbf{B}^{(1)}_{(2)}\right]$ will have possible non-zero structure constants $C^{l_{2}}_{i_{2}j_{2}}$. As we can observe, due to \eqref{be14} and the fact that $l_{3}=l_{1}+l_{2}$, $\grm_{1}=i_{2}+i_{4}$ we come to the conclusion that $C^{l_{2}}_{i_{2}j_{2}}=0$. This implies
\begin{align}
\left[\mathbf{B}^{(1)}_{(2)},\mathbf{B}^{(1)}_{(2)}\right]=0.\label{be18}
\end{align}
As a result, $\mathbf{B}^{(0)}_{(2)}$ is a (2)-level solvable algebra, and for simplicity we could define, $\mathbf{\bar{g}}^{(0)}=\mathbf{B}^{(0)}_{(2)}$, $\mathbf{\bar{g}}^{(1)}=\left[\mathbf{B}^{(0)}_{(2)},\mathbf{B}^{(0)}_{(2)}\right]$, $\mathbf{\bar{g}}^{(2)}=\left[\mathbf{B}^{(1)}_{(2)},\mathbf{B}^{(1)}_{(2)}\right]$ so that to the solvability chain reads 
\begin{align}
\left\{0\right\}\equiv\mathbf{\bar{g}}^{(2)}\subset\mathbf{\bar{g}}^{(1)}\subset\mathbf{\bar{g}}^{(0)}.\label{be19}
\end{align}
Now in the previous section we have already show that in the (2)-level solvable algebra we can use all the generators. Therefore, after the ``first integration method'' applied in $\mathbf{B}^{(1)}_{(2)}$, we are left with the generators $\mathbf{B}^{(0)}_{(1)}=\mathbf{B}^{(0)}_{(2)}-\mathbf{B}^{(1)}_{(2)}$, which will form an Abelian algebra, leading to the following chain
\begin{align}
\mathbf{B}^{(2)}_{(3)}\xrightarrow{pr}\mathbf{B}^{(1)}_{(2)}\xrightarrow{pr}\mathbf{B}^{(0)}_{(1)}.\label{be20}
\end{align}
In order to avoid being tedious, we will not study any higher level solvable algebra. The procedure can be iterated (n)-times thus, this concludes the proof.

All the above information can be gathered into the basic theorem of this work.

\textbf{\underline{Theorem 5}:} \textit{\textbf{Let us consider a system of k ordinary differential equations of maximal order m, which admits an (n)-level solvable algebra of order r. Furthermore, let us call the dimension of the system to be $\mathbf{N=m\,c_{(m)}+(m-1)\,c_{(m-1)}+...+2\,c_{(2)}}+1\,c_{(1)}$ with $\mathbf{c_{(m)},c_{(m-1)},...,c_{(1)}}$ the number of equations of the corresponding order $\mathbf{m,(m-1),...,1}$, and $\mathbf{c_{(m)}+c_{(m-1)}+...+c_{(1)}=k}$, so that $\mathbf{N\geq r}$. The solution to this system can be found as, a solution of a system of ordinary differential equations of dimension $\mathbf{\bar{N}=N-r}$ and an (r)-number of quadratures provided via the ``first integration method'' acted, at least, (n)-times upon the chain of the prolonged and reduced ``cosets'' }} $\left(\mathbf{B^{(i)}_{(j)}=\mathbf{g}^{(i)}-\mathbf{g}^{(j)}}, \mathbf{g}^{(i)},\mathbf{g}^{(j)} \textbf{the derived algebras of level}\,(i), (j)\right)$ \textbf{correspondingly}
\begin{align}
\mathbf{B}^{(n-1)}_{(n)}\xrightarrow{pr}\mathbf{B}^{(n-2)}_{(n-1)}\xrightarrow{pr}....\xrightarrow{pr}\mathbf{B}^{(0)}_{(1)}.\label{be21}
\end{align}
\textbf{\textit{The procedure should start with $\mathbf{B}^{(n-1)}_{(n)}$, given that the elements of each ``coset" act transitive on some sub-space of the dependent variables, of dimension equal to the number of the corresponding ``coset". The symbol ``pr'' above the arrow implies that at each ``coset" except $\mathbf{B}^{(n-1)}_{(n)}$ we have to consider the properly reduced prolongation of it's generators.}}

The phrase ``at least'' refers to the fact that this is the minimum number of steps that we have to perform in order to use the full set of generators. For instance, since at every step the generators of the ``coset" form an Abelian algebra, we could either transform them all at once into their normal forms, as we considered above, or by repeated steps. Either way, the result will be the same, but in the second case with much more steps. We may note here the following situation: Assume one differential equation of dimension ($n$) which admits two symmetry generators forming an Abelian algebra. If we transform both generators into normal form at once, one along the dependent variable and the other along the independent, then the equation is reduced to dimension $(n-1)$ and not $(n-2)$. The reason for this is that the generator normal along the independent variable does not reduce the order, but rather transforms the equation into autonomous i.e. does not depend on the independent variable. Thus, in such kind of cases, it is more profitable to use more than the minimum number of steps.     

\section{Example}

It is instructive to present an example so as to make clear the statements of the final theorem. Let us consider the following system of differential equations
\begin{align}
&x^{''}(t)-\frac{(x^{'}(t))^{2}}{x(t)}=0,\label{ex1}\\
&y^{'''}(t)+y(t)-e^{-t}\frac{x^{'}(t)}{x(t)}=0,\label{ex2}
\end{align}
where ${(')}$ indicates first derivative with respect to $t$. For simplicity, we assign the symbols $x^{''}\rightarrow x_{2}, x^{'}\rightarrow x_{1}$ and so forth, while we also omit the time dependence $(t)$
\begin{align}
&x_{2}-\frac{x_{1}^{2}}{x}=0,\label{ex3}\\
&y_{3}+y-e^{-t}\frac{x_{1}}{x}=0.\label{ex4}
\end{align}
It can be shown that this system admits, the following Lie point symmetry generators
\begin{align}
\mathbf{Z}_{1}=x\partial_{x},\,\,\mathbf{Z}_{2}=x\partial_{x}+e^{-t}\partial_{y},\,\,\mathbf{Z}_{3}=\partial_{t}+x \ln{x}\,\partial_{x},\label{es1}
\end{align}
with their algebra being
\begin{align}
\left[\mathbf{Z}_{1},\mathbf{Z}_{2}\right]=0,\,\,\left[\mathbf{Z}_{2},\mathbf{Z}_{3}\right]=\mathbf{Z}_{2},\,\,\left[\mathbf{Z}_{3},\mathbf{Z}_{1}\right]=-\mathbf{Z}_{1}.\label{al1}
\end{align}
This is a $(2)-level$ solvable algebra with $g^{(1)}=\left\{\mathbf{Z}_{1},\mathbf{Z}_{2}\right\},\,\,g^{(2)}=\left\{0\right\}$. The two ``coset" of the chain to be followed are $B^{(1)}_{(2)}=\left\{\mathbf{Z}_{1},\mathbf{Z}_{2}\right\}$,\,\, $B^{(0)}_{(1)}=\left\{\mathbf{Z}_{3}\right\}$. According to the Theorem $5$, the system has dimension $N=3*1_{(3)}+2*1_{(2)}=5$, and since there are $r=3$ symmetry generators, we can reduce the dimension of the system to $\bar{N}=N-r=2$. A system of dimension $2$, can consist of either two first order ordinary differential equations, or one algebraic and one second order ordinary differential equation. We will present both cases.

Based on Theorem $5$ the procedure must start with $B^{(1)}_{(2)}$. The two generators of this ``coset" act transitive in the two dimensional plane $(x,y)$ thus we can transform them into normal form at once. Let us introduce the new coordinates $(s,u,w)$ such that $s$ is the independent variable and $\tilde{\mathbf{Z}}_{1}=\partial_{u},\,\tilde{\mathbf{Z}}_{2}=\partial_{w}$. The transformation reads
\begin{align}
s=t,\,\,u=\ln{x}-e^{t}y,\,\,w=e^{t}y,\label{mfgfk}
\end{align}
while in those variables, the third vector field becomes
\begin{align}
\tilde{\mathbf{Z}}_{3}=\partial_{s}+u\partial_{u}+w\partial_{w}.\label{mgjfk1}
\end{align}
To this end, the system of the equations in the transformed variables reads
\begin{align}
&w_{3}-u_{1}+2w_{1}-3w_{2}=0,\label{lsdksk}\\
&u_{2}+w_{2}=0.\label{sdsdl1}
\end{align}
The reduction takes place once we are restricted to the hypersurfaces $(u=constant, w=constant,s=\grt, p=u_{1}, q=w_{1})$.
\begin{align}
&q_{2}-p+2q-3q_{1}=0,\label{dfdfd}\\
&p_{1}+q_{1}=0,\label{sdsdl}\\
&u=\int{p[\grt(s)]ds}+c_{1},\label{msksl34}\\
&w=\int{q[\grt(s)]ds}+c_{2}.\label{msksl35}
\end{align}
In order to use the ``coset" $B^{(0)}_{(1)}$ we have to calculate the first prolongation of $\tilde{\mathbf{Z}}_{3}$, and then reduce it to the variables $(\grt,p,q)$. The first prolongation is 
\begin{align}
\tilde{\mathbf{Z}}^{(1)}_{3}=\partial_{s}+u\partial_{u}+w\partial_{w}+u_{1}\partial_{u_{1}}+w_{1}\partial_{w_{1}},\label{sdsdshjjj}
\end{align}
while the reduced generator will be calculate as follows
\begin{align}
&\mathbf{Y}_{3}=\tilde{\mathbf{Z}}^{(1)}_{3}(\grt)\partial_{\grt}+\tilde{\mathbf{Z}}^{(1)}_{3}(p)\partial_{p}+\tilde{\mathbf{Z}}^{(1)}_{3}(q)\partial_{q}\Rightarrow\nonumber\\
&\mathbf{Y}_{3}=\tilde{\mathbf{Z}}^{(1)}_{3}(s)\partial_{\grt}+\tilde{\mathbf{Z}}^{(1)}_{3}(u1)\partial_{p}+\tilde{\mathbf{Z}}^{(1)}_{3}(w1)\partial_{q}\Rightarrow\nonumber\\
&\mathbf{Y}_{3}=\partial_{\grt}+p\partial_{p}+q\partial_{q}.\label{mmfnfm}
\end{align}
From \eqref{mmfnfm} we recognize that the reduced generator \eqref{mgjfk1} has the same form as the transformed one. This is not always the case, it is just a characteristic of the special form that the first prolongation has.

Now we are ready to continue the procedure. The new coordinates that will transform $\mathbf{Y}_{3}$ into it's normal form will be $\grs,m,h$ where $\grs$ the independent variable and $\tilde{\mathbf{Y}}_{3}=\partial_{h}$.
\begin{align}
\mathbf{Y}_{3}\grs(\grt,p,q)=0,\,\mathbf{Y}_{3}m(\grt,p,q)=0,\,\mathbf{Y}_{3}h(\grt,p,q)=1.\label{mng1}
\end{align}
The possible solutions to this system of partial differential equations read
\begin{align}
\grs(\grt,p,q)=f_{1}(e^{-\grt}p,e^{-\grt}q),\,m(\grt,p,q)=f_{2}(e^{-\grt}p,e^{-\grt}q),\,h(\grt,p,q)=\grt+f_{3}(e^{-\grt}p,e^{-\grt}q),\label{mng2}
\end{align}
where $f_{1},f_{2},f_{3}$ some arbitrary functions of their arguments.

\subsection{One second order and one algebraic differential equation}

In order to end up with a system of one second order and one algebraic differential equation we have to set,
\begin{align}
\grs=e^{-\grt}p,\,m=e^{-\grt}q,\,h=\grt,\label{msdks1}
\end{align}
with the inverse being
\begin{align}
\grt=h,\,p=e^{h}\grs\,q=e^{h}m.\label{msdk2}
\end{align}
Let us also calculate $p_{1},q_{1},q_{2}$.
\begin{align}
p_{1}=\frac{e^{h}\left(1+h_{1}\grs\right)}{h_{1}},\,q_{1}=\frac{e^{h}\left(m_{1}+h_{1}m\right)}{h_{1}},\,q_{2}=\frac{e^{h}\left(h_{1}m_{2}+2h_{1}^{2}m_{1}+h_{1}^{3}m-h_{2}m_{1}\right)}{h_{1}^{3}}.\label{mshk3}
\end{align}
Use these relations into \eqref{dfdfd},\eqref{sdsdl} and solve with respect to $m_{2},h_{1}$ to get
\begin{align}
&m_{2}-\frac{\grs-m_{1}\left\{m^{3}h_{2}-(2+m_{1})\grs+3m^{2}h_{2}\grs+h_{2}\grs^{3}+m\left[\left(1+m_{1}\right)^{2}+3h_{2}\grs^{2}\right]\right\}}{(1+m_{1})(m+\grs)^{2}}=0,\label{mlmn1}\\
&h_{1}+\frac{1+m_{1}}{m+\grs}=0.\label{mlml2}
\end{align}
By the introduction of the hypersurface's $h=constant$ variables $\grs=k,\,v=h_{1}$ we are left with the system
\begin{align}
&m_{2}-\frac{k-m_{1}\left\{m^{3}v_{1}-(2+m_{1})k+3m^{2}v_{1}k+v_{1}k^{3}+m\left[\left(1+m_{1}\right)^{2}+3v_{1}k^{2}\right]\right\}}{(1+m_{1})(m+k)^{2}}=0,\label{mlmn15}\\
&v+\frac{1+m_{1}}{m+k}=0,\label{mlml26}\\
&h=\int{v[k(\grs)]d\grs}+c_{3}.\label{malka27}
\end{align}
Thus, indeed the original system has been reduced to a second order ordinary differential equation for $m$ and an algebraic one for $\gry$, \eqref{mlmn15},\eqref{mlml26}, alongside with three quadratures \eqref{malka27},\eqref{msksl34},\eqref{msksl35}, where $c_{1},c_{2},c_{3}$ are integration constants. Note that since one of the equations is algebraic, we can solve it already and using it in the second order one, to get
\begin{align}
m_{2}+\frac{(1+m_{1})^{2}\left[2m m_{1}+(-1+m_{1})k\right]}{(m+k)^{2}}=0.\label{fso}
\end{align}

\subsection{Two first order differential equations}

In order to end up with two first order differential equations, we have to use the other possible transformation
\begin{align}
\grs=e^{-\grt}q,\,m=e^{-\grt}p,\,h=\grt,\label{msdks1mm}
\end{align}
with the inverse being
\begin{align}
\grt=h,\,p=e^{h}m,\,q=e^{h}\grs.\label{msdk2fgh}
\end{align}
Let us also calculate $p_{1},q_{1},q_{2}$.
\begin{align}
p_{1}=\frac{e^{h}\left(m_{1}+h_{1}m\right)}{h_{1}},\,q_{1}=\frac{e^{h}\left(1+h_{1}\grs\right)}{h_{1}},\,q_{2}=\frac{e^{h}\left(2h_{1}^{2}-h_{2}+h_{1}^{2}\grs\right)}{h_{1}^{3}}.\label{mshk3ghg}
\end{align}
Use these relations into \eqref{dfdfd},\eqref{sdsdl} and solve with respect to $m_{1},h_{2}$ to get
\begin{align}
&m_{1}+h_{1}(m+\grs)+1=0,\label{mlmn1ssss}\\
&h_{2}+h_{1}^{2}(1+h_{1}m)=0.\label{mlml2ffff}
\end{align}
By the introduction of the reduction variables $\grs=k,\,v=h_{1}$ we are left with the system
\begin{align}
&m_{1}+v(m+k)+1=0,\label{mlmn1sssfgfds}\\
&v_{1}+v^{2}(1+v\,m)=0,\label{mlml2fddfdfff}\\
&h=\int{v[k(\grs)]d\grs}+c_{3}.\label{malka27hgg}
\end{align}
Thus, indeed the original system has been reduced to a system of two first order differential equations, \eqref{mlmn1sssfgfds},\eqref{mlml2fddfdfff}, alongside with three quadratures \eqref{malka27hgg},\eqref{msksl34},\eqref{msksl35}, where $c_{1},c_{2},c_{3}$ are integration constants.

Which of these two systems of equations is to our best interest to end up with? Our sole purpose is to solve the original system of differential equations. In both cases we have used all the symmetries that we had into our disposal. What is better, one second order, \eqref{fso}, or two first order differential equations \eqref{mlmn1sssfgfds}, \eqref{mlml2fddfdfff}? At this point what we could do is to search for new symmetries of the reduced equations in both cases. Even though, a system of first order differential equations admits a large number of symmetries, it is practically impossible to find even one of them, as it was argued in \cite{stephani_1990}. On the other hand, we can relatively easy find if a second order differential equations admits any symmetry. For instance, in this case,\eqref{fso} admits the symmetry $\grj=k\partial_{k}+m\partial_{m}$ which by use of the method reduces the equation to a first order. Thus, we eventually end up with one first order differential equation against a system of two. Therefore, we might say that perhaps it would be to our best interest to use as much symmetries as possible in order to reduce some of the differential equations into algebraic ones.

\section{Maximal solvable sub-algebra}

So far, we have proven that the generators belonging to some solvable sub-algebra could all be used, given that the proper order is followed, for the integration procedure. Thus, is trivial to say that from all the solvable sub-algebras, it is to our best interest to use the maximal, meaning the one with the maximum number of generators. This might not be unique, but all of them will lead to the same dimension of reduction. 

At this section we would like to address the following question:

\textbf{Is the maximal solvable sub-algebra the most profitable way to use the symmetry generators? In other words, is it possible to use alongside with the maximal, solvable sub-algebra, generators which don't belong to this sub-algebra, via the ``first integration method''?}

Let us recall the ``inheritance conditions" for a set of generators ${\mathbf{Z}_{i},\mathbf{Z}_{\bar{l}}}$ where $\mathbf{Z}_{i}$ the generators that we are going to transform into normal form, while $\mathbf{Z}_{\bar{l}}$ the inherited ones.
\begin{align}
&[\mathbf{Z}_{\bar{l}},\mathbf{Z}_{\bar{k}}]=C^{\bar{q}}_{\bar{l}\bar{k}}\mathbf{Z}_{\bar{q}}+C^{l}_{\bar{l}\bar{k}}\mathbf{Z}_{l},\label{sm1}\\
&[\mathbf{Z}_{i},\mathbf{Z}_{\bar{k}}]=C^{l}_{i\bar{k}}\mathbf{Z}_{l},\label{sm2}\\
&[\mathbf{Z}_{i},\mathbf{Z}_{j}]=0.\label{sm3}
\end{align}
Due to \eqref{sm2} we could also say that the generators $\mathbf{Z}_{\bar{k}}$ belong to the ideal of $\mathbf{Z}_{i}$. Let us further recall the form of commutator relations for a (2)-level solvable algebra where $\mathbf{Z}_{\grm_{0}}\in g^{(0)}$,\,$\mathbf{Z}_{\grm_{1}}\in g^{(1)}$,\,$\mathbf{Z}_{i_{1}}\in \mathbf{B}^{(0)}_{(1)}$.
\begin{align}
&[\mathbf{Z}_{i_{1}},\mathbf{Z}_{j_{1}}]=C^{\grl_{1}}_{i_{1}j_{1}}\mathbf{Z}_{\grl_{1}},\label{sm4}\\
&[\mathbf{Z}_{\grm_{1}},\mathbf{Z}_{j_{1}}]=C^{\grl_{1}}_{\grm_{1}j_{1}}\mathbf{Z}_{\grl_{1}},\label{sm5}\\
&[\mathbf{Z}_{\grm_{1}},\mathbf{Z}_{\grn_{1}}]=0.\label{sm6}
\end{align}
The correspondence between those set of three relations is achieved for $\mathbf{Z}_{\grm_{1}}\rightarrow \mathbf{Z}_{i}$,\,$\mathbf{Z}_{i_{1}}\rightarrow \mathbf{Z}_{\bar{k}}$. We can easily note that the ``inheritance conditions" does not imply that the algebra is solvable, while on the other hand a solvable algebra will satisfy the ``inheritance conditions" as we have already proven in the previous sections.

In order to answer the question posted previously, we will follow two paths, since as we have proven the order of using the generators in the integration procedure is important. Furthermore, the proof could be easily iterated in higher order solvable sub-algebras.

\subsection{$\mathbf{Z}_{\grm_{0}}$ has a larger ideal}

Let us assume that there exists an algebra with generators $\{\mathbf{Z}_{\grm_{0}},\mathbf{Z}_{I}\}$ where $\{\mathbf{Z}_{\grm_{0}}\}$ will correspond, by assumption, to the maximal solvable sub-algebra, while $\mathbf{Z}_{I}$ to the generators outside this sub-algebra. The dimensions will be called $D_{s}$ for the solvable sub-algebra and $D_{r}$ for the rest of the generators. The question to be answered:

\textbf{Is it possible to find a sub-set of generators $\mathbf{Z}_{A}$ of the set $\mathbf{Z}_{I}$ that satisfy the ``inheritance conditions" with all $\mathbf{Z}_{\grm_{1}}$ and $\mathbf{Z}_{i_{1}}$, such that after the use of $\mathbf{Z}_{\grm_{0}}\equiv\{\mathbf{Z}_{\grm_{1}},\mathbf{Z}_{i_{1}}\}$, the $\mathbf{Z}_{A}$ could also be used for the integration procedure?}

The ``inheritance conditions" for the full set $\{\mathbf{Z}_{\grm_{1}},\mathbf{Z}_{i_{1}},\mathbf{Z}_{A}\}$ will read:
\begin{align}
&[\mathbf{Z}_{i_{1}},\mathbf{Z}_{j_{1}}]=C^{\grl_{1}}_{i_{1}j_{1}}\mathbf{Z}_{\grl_{1}},\label{sm7}\\
&[\mathbf{Z}_{\grm_{1}},\mathbf{Z}_{j_{1}}]=C^{\grl_{1}}_{\grm_{1}j_{1}}\mathbf{Z}_{\grl_{1}},\label{sm8}\\
&[\mathbf{Z}_{\grm_{1}},\mathbf{Z}_{\grn_{1}}]=0,\label{sm9}\\
&[\mathbf{Z}_{A},\mathbf{Z}_{B}]=C^{M}_{AB}\mathbf{Z}_{M}+C^{l_{1}}_{AB}\mathbf{Z}_{l_{1}}+C^{\grl_{1}}_{AB}\mathbf{Z}_{\grl_{1}},\label{sm10}\\
&[\mathbf{Z}_{A},\mathbf{Z}_{j_{1}}]=C^{M}_{A j_{1}}\mathbf{Z}_{M}+C^{l_{1}}_{A j_{1}}\mathbf{Z}_{l_{1}}+C^{\grl_{1}}_{A j_{1}}\mathbf{Z}_{\grl_{1}},\label{sm11}\\
&[\mathbf{Z}_{\grm_{1}},\mathbf{Z}_{B}]=C^{\grl_{1}}_{\grm_{1}B}\mathbf{Z}_{\grl_{1}}.\label{sm12}
\end{align}
Note two things: First, the equations \eqref{sm7},\eqref{sm8},\eqref{sm9} do not involve any term $\mathbf{Z}_{A}$ since by assumption, they form a proper sub-algebra (the maximal) which is also solvable. Second, the set $\{\mathbf{Z}_{\grm_{1}},\mathbf{Z}_{i_{1}},\mathbf{Z}_{A}\}$ does not form a solvable algebra, as can be deduced from \eqref{sm10},\eqref{sm11}. We continue with the reduction of the system by the steps we already explained in previous sections, in which the following will hold $\mathbf{Z}_{\grm_{1}}\rightarrow \mathbf{Y}_{\grm_{1}}=0$, $\mathbf{Z}_{i_{1}}\rightarrow\mathbf{Y}_{i_{1}}$, $\mathbf{Z}_{A}\rightarrow\mathbf{Y}_{A}$, while the structure constants of the reduced generators have the same values with the ones of the original generators:
\begin{align}
&\eqref{sm7}\Rightarrow[\mathbf{Y}_{i_{1}},\mathbf{Y}_{j_{1}}]=0,\label{sm13}\\
&\eqref{sm10}\Rightarrow[\mathbf{Y}_{A},\mathbf{Y}_{B}]=C^{M}_{AB}\mathbf{Y}_{M}+C^{l_{1}}_{AB}\mathbf{Y}_{l_{1}},\label{sm14}\\
&\eqref{sm11}\Rightarrow[\mathbf{Y}_{A},\mathbf{Y}_{j_{1}}]=C^{M}_{A j_{1}}\mathbf{Y}_{M}+C^{l_{1}}_{A j_{1}}\mathbf{Y}_{l_{1}}.\label{sm15}
\end{align}
All the other equations are identically satisfied. As we have said, by assumption $\mathbf{Z}_{A}$ must belong also to the ideal of the $\mathbf{Z}_{i_{1}}$ at the reduced level, therefore this implies that 
\begin{align}
C^{M}_{A j_{1}}=0,\label{sm16}
\end{align}
which brings the previous equations into the exact form of the ``inheritance conditions"  
\begin{align}
&[\mathbf{Y}_{i_{1}},\mathbf{Y}_{j_{1}}]=0,\label{sm17}\\
&[\mathbf{Y}_{A},\mathbf{Y}_{B}]=C^{M}_{AB}\mathbf{Y}_{M}+C^{l_{1}}_{AB}\mathbf{Y}_{l_{1}},\label{sm18}\\
&[\mathbf{Y}_{A},\mathbf{Y}_{j_{1}}]=C^{l_{1}}_{A j_{1}}\mathbf{Y}_{l_{1}}.\label{sm19}
\end{align}
The second reduction implies $\mathbf{Y}_{i_{1}}\rightarrow \mathbf{X}_{i_{1}}=0$, $\mathbf{Y}_{A}\rightarrow\mathbf{X}_{A}$, and the inherited generators satisfy the algebra
\begin{align}
[\mathbf{X}_{A},\mathbf{X}_{B}]=C^{M}_{AB}\mathbf{X}_{M}.\label{sm20}
\end{align}
From those generators, eventually we will be able to use only an Abelian sub-algebra, let us call it $\mathbf{X}_{\bar{A}}$, $\bar{A}=1,...,\gra$, which implies
\begin{align}
[\mathbf{X}_{\tilde{A}},\mathbf{X}_{\tilde{B}}]=0\Leftrightarrow C^{M}_{\tilde{A}\tilde{B}}=0.\label{sm21}
\end{align}
Thus, we were able to use for the integration the $D_{s}$ generators of the solvable sub-algebra and $D_{a}<D_{r}$ of the rest of the generators, which gives a total of $D_{s}+D_{a}$. However, a contradiction appears. Let us rewrite the equations \eqref{sm7} to \eqref{sm12} for the sub-set $\{\mathbf{Z}_{\grm_{1}},\mathbf{Z}_{i_{1}},\mathbf{Z}_{\tilde{A}}\}$ by taking also into account the relations that come up along the way \eqref{sm16},\eqref{sm21}.
\begin{align}
&[\mathbf{Z}_{i_{1}},\mathbf{Z}_{j_{1}}]=C^{\grl_{1}}_{i_{1}j_{1}}\mathbf{Z}_{\grl_{1}},\label{sm22}\\
&[\mathbf{Z}_{\grm_{1}},\mathbf{Z}_{j_{1}}]=C^{\grl_{1}}_{\grm_{1}j_{1}}\mathbf{Z}_{\grl_{1}},\label{sm23}\\
&[\mathbf{Z}_{\grm_{1}},\mathbf{Z}_{\grn_{1}}]=0,\label{sm24}\\
&[\mathbf{Z}_{\tilde{A}},\mathbf{Z}_{\tilde{B}}]=C^{l_{1}}_{\tilde{A}\tilde{B}}\mathbf{Z}_{l_{1}}+C^{\grl_{1}}_{\tilde{A}\tilde{B}}\,\mathbf{Z}_{\grl_{1}},\label{sm25}\\
&[\mathbf{Z}_{\tilde{A}},\mathbf{Z}_{j_{1}}]=C^{l_{1}}_{\tilde{A} j_{1}}\mathbf{Z}_{l_{1}}+C^{\grl_{1}}_{\tilde{A} j_{1}}\mathbf{Z}_{\grl_{1}},\label{sm26}\\
&[\mathbf{Z}_{\grm_{1}},\mathbf{Z}_{\tilde{B}}]=C^{\grl_{1}}_{\grm_{1}\tilde{B}}\mathbf{Z}_{\grl_{1}}.\label{sm27}
\end{align}
It is not difficult to obtained the result, that these relations correspond to a (3)-level solvable sub-algebra. This would imply that we have found a solvable sub-algebra of dimension $D_{s}+D_{a}>D_{s}$ which contradicts our original assumption that ${\mathbf{Z}_{\grm_{0}}}$ represent the generators of the maximal solvable sub-algebra. Hence, the following corollary can be stated:

\textbf{\underline{Corollary 4}:} \textit{\textbf{Once the maximal solvable sub-algebra has been obtained, there exist no generators outside this sub-algebra, that would belong to the ideal of the maximal solvable sub-algebra at each step of integration, and satisfy the ``inheritance conditions".}}

\subsection{$\mathbf{Z}_{\grm_{0}}$ belongs to an ideal}

The question to be answered can be stated as follows:

\textbf{Is it possible to find a sub-set of generators $\mathbf{Z}_{A}$ of the set $\mathbf{Z}_{I}$ such that $\mathbf{Z}_{\grm_{0}}$, $\mathbf{Z}_{I_{0}}\subseteq\mathbf{Z}_{I}$ belong to it's ideal?}

The ``inheritance conditions" for $\mathbf{Z}_{\grm_{0}}$, $\mathbf{Z}_{I_{0}}$ will read
\begin{align}
&[\mathbf{Z}_{A},\mathbf{Z}_{B}]=0,\label{sm28}\\
&[\mathbf{Z}_{A},\mathbf{Z}_{J_{0}}]=C^{M}_{AJ_{0}}\mathbf{Z}_{M},\label{sm29}\\
&[\mathbf{Z}_{A},\mathbf{Z}_{\grm_{0}}]=C^{M}_{A \grm_{0}}\mathbf{Z}_{M},\label{sm30}\\
&[\mathbf{Z}_{\grm_{0}},\mathbf{Z}_{\grn_{0}}]=C^{\grl_{0}}_{\grm_{0}\grn_{0}}\mathbf{Z}_{\grl_{0}},\label{sm31}\\
&[\mathbf{Z}_{\grm_{0}},\mathbf{Z}_{J_{0}}]=C^{\grl_{0}}_{\grm_{0}J_{0}}\mathbf{Z}_{\grl_{0}}+C^{K_{0}}_{\grm_{0}J_{0}}\mathbf{Z}_{K_{0}}+C^{M}_{\grm_{0}J_{0}}\mathbf{Z}_{M},\label{sm32}\\
&[\mathbf{Z}_{I_{0}},\mathbf{Z}_{J_{0}}]=C^{\grl_{0}}_{I_{0}J_{0}}\mathbf{Z}_{\grl_{0}}+C^{K_{0}}_{I_{0}J_{0}}\mathbf{Z}_{K_{0}}+C^{M}_{I_{0}J_{0}}\mathbf{Z}_{M}.\label{sm33}
\end{align}
Once more, \eqref{sm31} holds by the assumption that the set $\{\mathbf{Z}_{\grm_{0}}\}$ forms a proper, solvable sub-algebra.  At the first step we use the generators $\mathbf{Z}_{A}$ for the integration procedure which read $D_{a}$ in number. Then, at the reduced level we could forget all the $\mathbf{Z}_{I_{0}}$ generators and use only the $\mathbf{Z}_{\grm_{0}}\rightarrow\mathbf{Y}_{\grm_{0}}$. The total number would be $D_{a}+D_{s}$. However, there is a contradiction once more. The subset $\{\mathbf{Z}_{A},\mathbf{Z}_{\grm_{0}}\}$ of $\{\mathbf{Z}_{A},\mathbf{Z}_{\grm_{0}},\mathbf{Z}_{I_{0}}\}$ based on \eqref{sm28},\eqref{sm30},\eqref{sm31} forms a (3)-level solvable sub-algebra of dimension $D_{a}+D_{s}>D_{s}$ with the $g^{(2)}=\{\mathbf{Z}_{A}\}$.   This contradicts the assumption that we have started with the maximal solvable sub-algebra. The following corollaries are deduced based on the previous results:

\textbf{\underline{Corollary 5}:} \textit{\textbf{Once the maximal solvable sub-algebra has been obtained, there exist no generators outside this sub-algebra to which, the maximal solvable sub-algebra would belong to their ideal and satisfy the ``inheritance conditions".}}

\textbf{\underline{Corollary 6}:} \textit{\textbf{The set of generators that will eventually be used for the integration of a given system of ordinary differential equations, will form a solvable sub-algebra.}}

To end this section, the following theorem is presented whose proof is already obtained from the previous steps and Corollaries.

\textbf{\underline{Theorem 6}:} \textit{\textbf{The use of the maximal solvable sub-algebra corresponds to the most profitable way of using the symmetry generators of a system of ordinary differential equations, in order to integrate it, based on the ``first integration method''.}}

\subsection{``Optimum'' maximal solvable sub-algebra}

At this section we would like to note something about the dimension of the maximal, solvable sub-algebra.

Suppose that we are given an algebra which is not solvable. Finding the maximal, solvable sub-algebra, depends on the representation of the structure constants in the given equivalence class. To make this clear, consider the following example.
\begin{align}
[\mathbf{X}_{1},\mathbf{X}_{2}]=-\mathbf{X}_{3},\,\,[\mathbf{X}_{1},\mathbf{X}_{3}]=-\mathbf{X}_{2},\,\,[\mathbf{X}_{2},\mathbf{X}_{3}]=\mathbf{X}_{1}.\label{mm1}
\end{align}  
This algebra is not solvable since $g^{(1)}=\{\mathbf{X}_{1},\mathbf{X}_{2},\mathbf{X}_{3}\}=g^{(0)}$. Based on the definition of solvability, the maximal solvable sub-algebra would be either one of these three $\{\mathbf{X}_{1}\}$,$\{\mathbf{X}_{2}\}$,$\{\mathbf{X}_{3}\}$, meaning that it's dimension would be $D_{s}=1$. However, by the following acceptable (meaning inside the same equivalent Class) change of basis, 
\begin{align}
\mathbf{X}_{1}=-\mathbf{Y}_{2},\,\,\mathbf{Z}_{2}=\mathbf{Y}_{3},\,\,\mathbf{Z}_{3}=\mathbf{Y}_{1}+\mathbf{Y}_{3},\label{mm2}
\end{align} 
it is not hard to see that this Class possesses a two dimensional sub-algebra
\begin{align}
[\mathbf{Y}_{1},\mathbf{Y}_{2}]=\mathbf{Y}_{1},\,\,[\mathbf{Y}_{1},\mathbf{Y}_{3}]=\mathbf{Y}_{2},\,\,[\mathbf{Y}_{2},\mathbf{Y}_{3}]=\mathbf{Y}_{1}+\mathbf{Y}_{3}.\label{mm3}
\end{align}
It is true that each two dimensional sub-algebra is solvable, therefore it is not hard to see that $\{\mathbf{Y}_{1},\mathbf{Y}_{2}\}$ forms the maximal solvable sub-algebra and has dimension $D_{s}=2$.

To conclude, in everything that we have proven so far, when we refer to the maximal, solvable sub-algebra is this ``optimum'' as we name it, maximal, solvable sub-algebra.

The way to really track it down is the following: Find the basis of vector fields in which the maximal sub-algebra is manifest. Check if this sub-algebra is also solvable. If not, for this sub-algebra we find it's maximal sub-algebra and we repeat the procedure until we find the first solvable sub-algebra. This would also correspond to the ``optimum'', maximal, solvable sub-algebra.

\section{Discussion}

In the present work we have focused on, an already known, specific integration procedure of ordinary differential equations. This procedure (as well as the others) relies on the existence of symmetry transformations of the equations; meaning transformations of independent and dependent variables that map solutions into solutions. These transformations can be best understood in terms of their infinitesimal generators which are realized as vector fields on the space of independent and dependent variables; they form a Lie algebra. This specific procedure we call ``first integration method'' and concerns finding new variables in which the transformed vector fields acquire their normal form along the new dependent variables.

In Theorem 1, we concentrate on the application of the ``first integration method'' upon an Abelian sub-algebra with more than one symmetry generators: we prove the properties that a generator $\mathbf{Z}_{\bar{l}}$ should have so that it is inherited as symmetry in the reduced system of equations. There are four conditions: the first expresses the fact that $\mathbf{Z}_{\bar{l}}$ should belong to the ideal of $\mathbf{Z}_{i}$. The second is concerned with the form of the reduced operator $\hat{\mathbf{A}}$ that associated to the reduced system of equations. The third and fourth are related to the form and the algebra that the properly reduced generators of $\mathbf{Z}_{\bar{l}}$, namely $\mathbf{Y}_{\bar{l}}$, satisfy. As a corollary, we deduce that the set $\{\mathbf{Z}_{i},\mathbf{Z}_{\bar{l}}\}$ forms a sub-algebra of the original algebra. 

The combination of solvability, a property of the generators algebra, and the first integration method leads to Theorem 5. This states that we can use all the generators of some solvable sub-algebra if we follow a specific order of application of the method: start from the last ``coset" all the way to the first one. At each step we must use the proper reduced and prolonged generators of the corresponding ``coset". In the case that the dimension of the system ($N$, see Theorem 5) is equal or less than the dimension of the algebra, the solution to the system of equations can be found in terms of quadratures only.

Another important result stemming out of this work is that the ``optimum'', maximal, solvable sub-algebra is the best we can get (see section VI). That is to say, there is no other way that we can use more generators for the integration procedure (via the ``first integration method'').

To conclude, we believe that this work provides the missing (in the literature) link between the combination of the property of solvability with the ``first integration method''. As we have stated, this combination has been established in the other two integration methods. Lastly, we believe that these theorems will prove to be useful for any application of this integration method to some system of ordinary differential equations, possible even outside the realm of physics.

\begin{appendices}

\section{Transitivity of the Abelian sub-algebra}

Let us assume the existence of an Abelian algebra of dimension $m$, with generators $\mathbf{Z}_{i}$ which acts on a space of dimension $n\geq m$. The generators are linearly independent, that is the equation
\begin{align}
a^{(i)}\mathbf{Z}_{i}=0,\label{ea1}
\end{align}
where $a^{(i)}$ are some constants admits the trivial and unique solution
\begin{align}
a^{(i)}=0,\forall i=1,...,m.\label{ea2}
\end{align}
The transitivity properties of the algebra can be found by whether there exist non-trivial solutions of the equation 
\begin{align}
a^{(i)}(x^{\grm})\mathbf{Z}_{i}=0,\label{ea3}
\end{align}
where $a^{(i)}(x^{\grm})$ are some functions, and $i=1,...,m$, $\grm=1,...,n$. For the purpose of a better understanding, let us consider the following example: Consider some space of dimension $n=2$ with coordinates $(t,s)$ and a two dimensional $m=2$ Abelian algebra, with realization 
\begin{align}
\mathbf{Z}_{1}=\partial_{s},\hspace{0.2cm} \mathbf{Z}_{2}=t\partial_{s}.\label{ea4} 
\end{align}
It is easy to verify that those generators form an Abelian algebra.
\begin{align}
\left[\mathbf{Z}_{1},\mathbf{Z}_{2}\right]=0.\label{ea5}
\end{align}
Furthermore, let us prove that those two are also linearly independent.
\begin{align}
&a^{(1)}\mathbf{Z}_{1}+a^{(2)}\mathbf{Z}_{2}=0\Rightarrow\nonumber\\
&(a^{(1)}+a^{(2)}t)\partial_{s}=0\Rightarrow\nonumber\\
&a^{(1)}+a^{(2)}t=0\Rightarrow\nonumber\\
&a^{(1)}=0,\hspace{0.2cm} a^{(2)}=0.\label{ea6}
\end{align}
Although, this group acts intransitive since there exist at least one, non-trivial solution of the equation \eqref{ea3}.
\begin{align}
&a^{(1)}(t,s)\mathbf{Z}_{1}+a^{(2)}(t,s)\mathbf{Z}_{2}=0\Rightarrow\nonumber\\
&(a^{(1)}(t,s)+a^{(2)}(t,s)t)\partial_{s}=0\Rightarrow\nonumber\\
&a^{(1)}(t,s)+a^{(2)}(t,s)t=0\Rightarrow\nonumber\\
&a^{(1)}(t,s)=-t\,a^{(2)}(t,s).\label{ea7}
\end{align}
The basic effect of intransitiveness, is that we cannot find new coordinates $x=x(t,s), y=y(t,s)$ such that the generators would acquire the form
\begin{align}
\mathbf{\tilde{Z}}_{1}=\partial_{x},\hspace{0.2cm} \mathbf{\tilde{Z}}_{2}=\partial_{y}.\label{ea8}
\end{align} 
As we have stated in Theorem 1 of the main text, we require additionally to the Abelian nature of the sub-algebra, that should also act transitive in a subspace of dimension equal to the dimension of the sub-algebra. That is equivalent to the existence of the trivial solution for the equation \eqref{ea3}. 

\end{appendices}

\section*{Acknowledgements}

\begin{figure}[h!]
\centering
  \begin{subfigure}[h]{0.265\linewidth}
    \includegraphics[width=\linewidth]{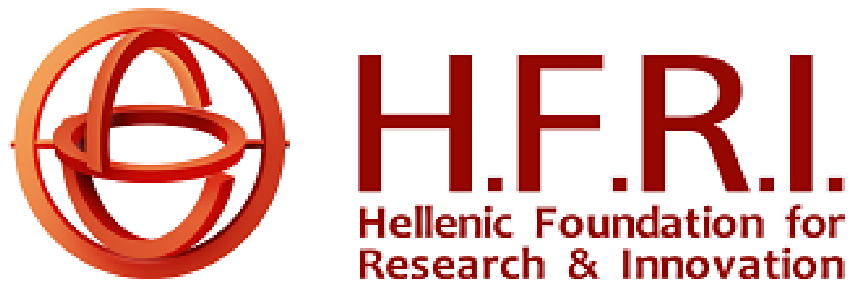}
  \end{subfigure}
  \begin{subfigure}[h]{0.2\linewidth}
    \includegraphics[width=\linewidth]{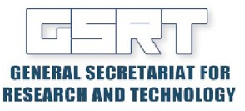}
  \end{subfigure}
\end{figure}
The research work was supported by the Hellenic Foundation for Research and Innovation (HFRI) and the General Secretariat for Research and Technology (GSRT), under the HFRI PhD Fellowship grant (GA.no.14501/2017).

\bibliographystyle{unsrt}

\begin{thebibliography}{10}

\bibitem{Sundermeyer:2014kha}
Kurt Sundermeyer.
\newblock {\em {Symmetries in fundamental physics}}, volume 176.
\newblock Springer, Cham, Switzerland, 2014.

\bibitem{Lie1}
S.~Lie.
\newblock {Bergundung einer Invariantentheorie der Beruhrungstransformationen}.
\newblock {\em Math. Ann.}, 8:215--288, 1874.

\bibitem{Lie2}
S.~Lie.
\newblock {Uber die Integration durch bestimmte Integrale von einer Klasse
  linear partieller Differentialgleichung}.
\newblock {\em Arch. for Math.}, 6:328--268, 1881.

\bibitem{Lie1888}
S.~Lie.
\newblock Classification und integration von gewöhnlichen
  differentialgleichungen zwischen xy, die eine gruppe von transformationen
  gestatten.
\newblock {\em Mathematische Annalen}, 32:213--281, 1888.

\bibitem{Lie1891}
Sophus Lie.
\newblock {\em Vorlesungen über Differentialgleichungen mit bekannten
  infinitesimalen Transformationen}.
\newblock Teubner, 1891.

\bibitem{Lie5}
S.~Lie.
\newblock Zur allgemeinen theorie der partiellen differentialgleichungen
  beliebeger ordnung.
\newblock {\em Berichte, Leipz.}, 47:53--128, 1895.

\bibitem{Lie6}
S.~Lie.
\newblock Die theorie der integralinvarianten ist ein korollar der theorie der
  differentialinvarianten.
\newblock {\em Berichte, Leipz.}, 49:342--357, 1897.

\bibitem{Ovs58}
L.~V. Ovsyannikov.
\newblock Groups and invariant-group solutions of differential equations.
\newblock {\em Dokl. Akad. Nauk USSR}, 118:439--442, 1958.

\bibitem{Ovs59}
L.~V. Ovsyannikov.
\newblock Group properties of the nonlinear heat conduction equation.
\newblock {\em Dokl. Akad. Nauk USSR}, 125:492--495, 1959.

\bibitem{Ovs60}
L.~V. Ovsyannikov.
\newblock {\em Group Properties of Differential Equations}.
\newblock Novosobirsk, 1962.

\bibitem{Ovs61}
L.~V. Ovsyannikov.
\newblock {\em Group Analysis of Differential Equations}.
\newblock Academic Press, 1982.

\bibitem{B0}
G.~W. Bluman and J.~D. Cole.
\newblock The general similarity solutions of the heat equation.
\newblock {\em J. Math. Mech.}, 18:1025--1042, 1969.

\bibitem{B1}
G.~W. Bluman.
\newblock Similarity solutions of the one-dimensional fokker-planck equation.
\newblock {\em Int. J. Nonlin. Mech.}, 6:143--153, 1971.

\bibitem{B2}
G.~W. Bluman.
\newblock Applications of the general similarity solution of the heat equation
  to boundary value problems.
\newblock {\em Quart. Appl. Math.}, 31:403--415, 1974a.

\bibitem{B3}
G.~W. Bluman.
\newblock Use of group methods for relating linear and nonlinear partial
  differential equations.
\newblock {\em Proceedings of Symposium on Symmetry, Similarity and Group
  Theoretic Methods in Mechanics}, pages 203--218, 1974b.

\bibitem{B4}
S.~Kumei and G.~W. Bluman.
\newblock When nonlinear differential equations are equivalent to linear
  differential equaitons.
\newblock {\em SIAM J. Appl. Math.}, 42:1157--1173, 1982.

\bibitem{B5}
F.~Schwarz.
\newblock Symmetries of the two dimensional korteweg-de vries equation.
\newblock {\em J. Phys. Soc. Japan}, 51:2387--2388, 1982.

\bibitem{Doyle}
J.~Doyle and M.~J. Englefield.
\newblock {Similarity Solutions of a Generalized Burgers Equation}.
\newblock {\em IMA Journal of Applied Mathematics}, 44(2):145--153, 1990.

\bibitem{Noether1918}
E.~Noether.
\newblock Invariante variationsprobleme.
\newblock {\em Nachrichten von der Gesellschaft der Wissenschaften zu
  Göttingen, Mathematisch-Physikalische Klasse}, 1918:235--257, 1918.

\bibitem{Tsamparlis:2011wf}
Michael Tsamparlis and Andronikos Paliathanasis.
\newblock {Lie and Noether symmetries of geodesic equations and collineations}.
\newblock {\em Gen. Rel. Grav.}, 42:2957--2980, 2010.

\bibitem{Paliathanasis:2017ocj}
Andronikos Paliathanasis.
\newblock {Dust fluid component from Lie symmetries in Scalar field Cosmology}.
\newblock {\em Mod. Phys. Lett.}, A32(37):1750206, 2017.

\bibitem{Paliathanasis:2011jq}
Andronikos Paliathanasis, Michael Tsamparlis, and Spyros Basilakos.
\newblock {Constraints and analytical solutions of $f(R)$ theories of gravity
  using Noether symmetries}.
\newblock {\em Phys. Rev.}, D84:123514, 2011.

\bibitem{Christodoulakis:2006vi}
T.~Christodoulakis and Petros~A. Terzis.
\newblock {The General solution of Bianchi type III vacuum cosmology}.
\newblock {\em Class. Quant. Grav.}, 24:875--887, 2007.

\bibitem{Terzis:2008ev}
Petros~A. Terzis and T.~Christodoulakis.
\newblock {The General Solution of Bianchi Type VII(h) Vacuum Cosmology}.
\newblock {\em Gen. Rel. Grav.}, 41:469--495, 2009.

\bibitem{Terzis:2010dk}
Petros~A. Terzis and T.~Christodoulakis.
\newblock {Lie algebra automorphisms as Lie point symmetries and the solution
  space for Bianchi Type I, II, IV, V vacuum geometries}.
\newblock {\em Class. Quant. Grav.}, 29:235007, 2012.

\bibitem{Pailas:2018tzy}
T.~Pailas, Petros~A. Terzis, and T.~Christodoulakis.
\newblock {The solution space of the Einstein’s vacuum field equations for
  the case of five-dimensional Bianchi Type I (Type 4A1)}.
\newblock {\em Class. Quant. Grav.}, 35(14):145003, 2018.

\bibitem{Christodoulakis:2013xha}
T.~Christodoulakis, N.~Dimakis, and Petros~A. Terzis.
\newblock {Lie point and variational symmetries in minisuperspace Einstein
  gravity}.
\newblock {\em J. Phys.}, A47:095202, 2014.

\bibitem{Bac1}
A.~V. Backlund.
\newblock Einiges uber curven und flachen transformationen.
\newblock {\em Afdelningen for Mathematik och Naturetenskap}, II:1--2, 1873.

\bibitem{Bac2}
A.~V. Backlund.
\newblock Ueber flachentransformationen.
\newblock {\em Math. Ann.}, IX:297--320, 1876.

\bibitem{Bac3}
A.~V. Backlund.
\newblock Zur theorie der partiellen differentialgleichungen erster ordnung.
\newblock {\em Math. Ann.}, XVII:285--328, 1880.

\bibitem{Bac4}
A.~V. Backlund.
\newblock Zur theorie der flachentransformationen.
\newblock {\em Math. Ann.}, XIX:387--422, 1882.

\bibitem{LE}
S.~Lie and F.~Engel.
\newblock Theorie der transformationsgruppen ii.
\newblock {\em Teubner Leipzig}, 1890.

\bibitem{Dimakis:2015rba}
N.~Dimakis, Petros~A. Terzis, and T.~Christodoulakis.
\newblock {Contact symmetries of constrained systems and the associated
  integrals of motion}.
\newblock {\em J. Phys. Conf. Ser.}, 633(1):012040, 2015.

\bibitem{G1}
Robert Geroch.
\newblock A method for generating solutions of einstein's equations.
\newblock {\em Journal of Mathematical Physics}, 12(6):918--924, 1971.

\bibitem{G2}
Robert Geroch.
\newblock A method for generating new solutions of einstein's equation. ii.
\newblock {\em Journal of Mathematical Physics}, 13(3):394--404, 1972.

\bibitem{Ibra1}
R.~L. Anderson and N.~H. Ibragimov.
\newblock {\em Lie-Backlund transformations in applications}.
\newblock Society for Industrial and Applied Mathematics, 1979.

\bibitem{Ibra2}
N.~H. Ibragimov.
\newblock {\em Transformation Groups Applied to Mathematical Physics}.
\newblock Cambridge University Press, 1985.

\bibitem{olver2000applications}
P.J. Olver.
\newblock {\em Applications of Lie Groups to Differential Equations}.
\newblock Applications of Lie Groups to Differential Equations. Springer New
  York, 2000.

\bibitem{Olv1}
P.~J. Olver.
\newblock Evolution equations possessing infinitely many symmetries.
\newblock {\em J. Math. Phys.}, 18:1212--1215, 1977.

\bibitem{Mae}
SHIGERU MAEDA.
\newblock {The Similarity Method for Difference Equations}.
\newblock {\em IMA Journal of Applied Mathematics}, 38(2):129--134, 05 1987.

\bibitem{LEVI1991335}
D.~Levi and P.~Winternitz.
\newblock Continuous symmetries of discrete equations.
\newblock {\em Physics Letters A}, 152(7):335 -- 338, 1991.

\bibitem{Levi_2010}
D.~Levi, P.~Winternitz, and R.~I. Yamilov.
\newblock Lie point symmetries of differential-difference equations.
\newblock {\em Journal of Physics A: Mathematical and Theoretical},
  43(29):292002, jun 2010.

\bibitem{QUISPEL1992379}
G.~R.~W. Quispel, H.~W. Capel, and R.~Sahadevan.
\newblock Continuous symmetries of differential-difference equations: the
  kac-van moerbeke equation and painlevé reduction.
\newblock {\em Physics Letters A}, 170(5):379 -- 383, 1992.

\bibitem{stephani_1990}
Hans Stephani.
\newblock {\em Differential Equations: Their Solution Using Symmetries}.
\newblock Cambridge University Press, 1990.

\bibitem{bluman2013symmetries}
G.W. Bluman and S.~Kumei.
\newblock {\em Symmetries and Differential Equations}.
\newblock Applied Mathematical Sciences. Springer New York, 2013.

\bibitem{Gilmore:102082}
Robert Gilmore.
\newblock {\em {Lie groups, Lie algebras, and some of their applications}}.
\newblock Wiley-Interscience, New York, NY, 1974.
\newblock Also a reprint ed.: Mineola, NY, Dover, 2005.

\end{thebibliography}

\end{document}